\documentclass[acmlarge,nonacm]{acmart}
\setcopyright{none}
\usepackage{amsmath,mathtools}
\usepackage{bbm} 
\usepackage{tikz-cd}
\usepackage{booktabs}
\usepackage{listings}

\newcommand{\Cxx}{C++}
\newcommand{\Lean}{Lean}

\newtheorem{theorem}{Theorem}[section]
\newtheorem{lemma}[theorem]{Lemma}
\newtheorem{proposition}[theorem]{Proposition}

\theoremstyle{definition}
\newtheorem{definition}[theorem]{Definition}
\theoremstyle{remark}
\newtheorem{remark}[theorem]{Remark}

\title{Semantics-Preserving Compilers via Operadic Algebras and Factorization Algebras }
\author{Yuwen Gu}
\affiliation{%
  \department{Shanghai Center for Mathematical Sciences}
  \institution{Fudan University}
  \city{Shanghai}
  \country{China}}
\email{ywgu24@m.fudan.edu.cn}
\authorsaddresses{Shanghai Center for Mathematical Sciences, Fudan University, Shanghai, China.
  \href{mailto:ywgu24@m.fudan.edu.cn}{\nolinkurl{ywgu24@m.fudan.edu.cn}}}

\begin{abstract}
We investigate an operadic approach to semantics-preserving translation from a well-defined subset of \Cxx{} to \Lean{}.
The proposed development has three components.
First, it organizes program constructors and their interpretations using colored operads and their algebras.
We give sufficient conditions for a correspondence on fundamental types and basic operations to induce a semantically faithful compiler, and prove that such a compiler extends through finite products and recursion.
Second, we use a factorization algebra to describe how the semantics of a program is assembled from its function fragments.
The local data include the embedding of each fragment algebra into the program's semantic algebra. We prove that, for a cover by function fragments, the images of these embeddings generate that algebra.
Third, we propose a compiler architecture based on the factorization algebra of program fragments and implement it for translation from \Cxx{} to \Lean{}. The compiler records each function as a fragment, translates each fragment using the correspondence on types and operations, and assembles the resulting definitions according to their call dependencies.
We compile and verify the Milnor multiplication program used by Lin, Wang, and Xu in their computations for the Kervaire invariant one problem in dimension~126. After repairing an uninitialized read, we verify that the translated program computes the product specified by the formalized Milnor coproduct for every ordered basis pair of combined degree at most~180.
As another application, we rewrite a \Lean{}-verified algorithm for ML-KEM polynomial multiplication in the \Lean{} fragment accepted by our compiler.
We prove in \Lean{} that the two formulations agree on all canonical inputs and compile the rewritten definition to \Cxx{}.
On an Apple M4, the generated \Cxx{} code ran 6.1 times as fast as the native \Lean{} version of the same definition.
\end{abstract}

\begin{document}
    
\maketitle

\section{Introduction}
\label{sec:introduction}
Computer programming has long played an important role in mathematical research. Over the past decades, numerous mathematical computing systems, such as \texttt{Mathematica}, \texttt{Maple}, \texttt{SageMath}, and \texttt{MATLAB}, have been developed to support symbolic and numerical computations. More recently, Lin, Wang, and Xu employed computer-assisted computations, including software written in C++, to resolve the final open case of the Kervaire invariant one problem, a long-standing problem in algebraic topology.

Despite these successes, the reliability of mathematical software remains an important concern. First, even when the underlying mathematical algorithms are correct, their implementations may contain bugs or deviate from the intended mathematical formulations, potentially leading to incorrect results. Second, mathematical software systems can be large and complex, making it difficult to understand and verify their behavior. Third, many existing computational workflows are not designed with formal verification in mind, making it challenging to establish rigorous connections between mathematical formulations and their implementations.

These considerations motivate us to develop a compiler that translates programs written in a well-defined subset of \Cxx{} into \Lean{}, a proof assistant and programming language based on dependent type theory, and to establish the semantic correctness of the compiler mathematically. To achieve this, we need to formulate the semantics of programming languages and semantics-preserving translations in mathematical terms.

The use of algebraic methods to study compiler correctness has a long history. Early work by Thatcher, Wagner, and Wright~\cite{thatcher-wagner-wright-compiler} established an algebraic approach to compiler construction and correctness proofs, while Mosses~\cite{mosses-constructive} investigated how correct compilers could be constructed from denotational semantics. These developments demonstrated the potential of algebraic methods for compiler construction. However, turning abstract semantic descriptions into practical compilation methods presented further difficulties. In particular, Mosses's Semantics Implementation System (SIS) was originally intended to generate correct-by-construction interpreters and compilers, but attempts to describe larger programming languages exposed significant limitations in the modularity of semantic descriptions~\cite{mosses-sis}.

These developments motivate two questions that are central to our work. First, given the algebraic semantics of two programming languages, how can we construct a compiler that preserves their semantics, particularly when their types, operations, and treatment of undefined behavior differ? Second, how can such a mathematical construction be organized using local program information, rather than requiring the entire program to be expanded before translation?

Roughly speaking, a computer program is a machine that takes inputs and produces outputs after performing a series of operations. Each elementary operation can be viewed as taking a finite number of inputs and producing an output. More importantly, complex programs are constructed by composing elementary operations: the output of one operation can be substituted into an input of another. Therefore, a mathematical model of programming languages should describe not only individual operations, but also how these operations can be composed.

An operad is a mathematical structure designed to describe precisely this kind of composition. It consists of collections of operations with specified numbers of inputs, together with composition rules that allow operations to be substituted into one another. For programming languages, we must also take types into account, since not every output can be substituted into every input. This leads naturally to colored operads, where each input and output is assigned a color representing its type.

The advantage of this approach is that the mathematical structure directly captures typed operations and their compositions in a unified algebraic framework. Moreover, operads and their algebras provide a natural setting for describing both the semantics of programming languages and the translations between them through structure-preserving morphisms.

To address the first question, we construct semantic algebras for programming languages using colored operads and their algebras, explicitly accounting for inputs on which a basic operation is undefined. We then formulate compilers as morphisms of operadic algebras using extension by zero. We establish sufficient conditions on the correspondences of generators, their defining relations, and semantic interpretations, together with suitable injectivity and closure conditions, under which a compiler can be constructed. We prove that the semantic representations of a source program and its translation are isomorphic. Furthermore, we show that the correctness of such a compiler is preserved under finite-product and recursion constructions, allowing the translation to be extended from fundamental types and operations to more complicated structures.

To address the second question, we introduce factorization algebras to organize the semantic information of local program fragments. We construct a factorization algebra from the semantic algebras of these fragments and establish a local-to-global reconstruction theorem. For an open cover by function fragments, extension by zero realizes the fragment algebras as subalgebras of the program's semantic algebra, and their join is the whole algebra. This provides a mathematical basis for organizing compilation through local information and its assembly, without repeatedly expanding function definitions at their call sites.

We apply these constructions to the public Milnor multiplication program in \texttt{SSeqCpp}, archived by Lin, Wang, and Xu with their computations for the Kervaire invariant one problem in dimension~126~\cite{lin-wang-xu-sseq}. We replace the pointer used to read the exponents of a monomial with indexed access and repair a boundary test in the original matrix search, \texttt{MulMilnorV3}, that permitted a read from an uninitialized cell. For every ordered basis pair of combined degree at most~180, we verify that the program computes the product specified by the Milnor coproduct formalized in \texttt{Kervaire126}~\cite{gair-kervaire126}.

We also use the compiler in the other direction.
For polynomial multiplication in ML-KEM, we start from a \Lean{}-verified algorithm based on the number-theoretic transform.
We express its transform layers using finite products and bounded recursion so that the algorithm lies in the \Lean{} fragment accepted by our compiler.
We prove in \Lean{} that the rewritten definition agrees with the sequential algorithm on every pair of canonical inputs and computes their product in \(\mathbb{F}_{3329}[X]/(X^{256}+1)\).
On an Apple M4, the \Cxx{} code generated from this definition ran 6.1 times as fast as the same definition compiled natively with \Lean{}.

\paragraph*{Organization}
Section~\ref{sec:semantic-algebra} constructs the semantic algebras of programs and programming languages and defines what it means for a language to represent a program faithfully.
Section~\ref{sec:compiler} gives sufficient conditions for constructing a semantically faithful compiler, and Section~\ref{subsec:compiler-extension} shows how to extend it to types and operations built by finite products and recursion.
Section~\ref{sec:local-to-global} uses a factorization algebra to recover the semantics of a program from its function fragments and the calls between them.
Section~\ref{sec:implementation} describes how the compiler records and translates each fragment, including its variable lifetimes, control flow, and parameter bindings.
Section~\ref{sec:compiler-applications} illustrates the two directions of translation with ML-KEM polynomial multiplication (\Lean{} to \Cxx{}) and Milnor multiplication (\Cxx{} to \Lean{}).
Section~\ref{sec:related-work-conclusion} relates our construction to earlier work on program semantics and compiler correctness and concludes the paper.

\paragraph*{Disclosure of AI use}
We independently developed the mathematical formulation presented in this paper, including its definitions, results, and proofs.
We designed the compiler architecture. This included deciding what information to extract from \Cxx{} functions, how to store the resulting function fragments in SQLite, how to translate their records using a dictionary of \Cxx{} types and operations and their \Lean{} counterparts, and how to assemble the translated functions into a complete \Lean{} program.

We used ChatGPT to help locate references, summarize prior work, refine the wording of the manuscript, and review the mathematical arguments. Following our design, ChatGPT also wrote implementation code, including code used for verification in \Lean{}, and ran tests for the examples reported in the paper.
We reviewed the AI-generated code. In the extensive tests reported in Section~\ref{sec:compiler-applications}, the final implementations produced the expected output for every input tested.

\section{Semantic Algebra}
\label{sec:semantic-algebra}
Our aim in this section is to find an algebraic framework for describing both the syntax and the semantics of a program. The framework is independent of the choice of programming language. We illustrate semantic restriction and semantic faithfulness through the example of a program represented within a language. This example motivates our definition of a compiler in the next section.
\subsection{Operads and Colored Operads}
\label{subsec:colored-operads}
Operads provide an algebraic framework for encoding operations with multiple inputs and a single output, together with their composition laws. Originally developed in algebraic topology through the work of Boardman–Vogt and May on iterated loop spaces and homotopy-invariant algebraic structures~\cite{boardman-vogt-homotopy,may-iterated-loop-spaces}, operads have since become a fundamental tool for describing compositional algebraic structures. The following definition of operads is adapted from May~\cite[Section~1]{may-iterated-loop-spaces} with some modifications to suit our context.
\begin{definition}
A \emph{operad} $\mathcal O$ consists of
a family of sets
\[
    \{\mathcal O(n)\}_{n\geq 0},
\]
together with the following data:

\begin{enumerate}
    \item \textbf{Composition.}
    For every $k\geq 0$ and $n_1,\ldots,n_k\geq 0$,
    there is a morphism
    \[
    \gamma:
    \mathcal O(k)\times
    \mathcal O(n_1)\times\cdots\times
    \mathcal O(n_k)
    \longrightarrow
    \mathcal O(n_1+\cdots+n_k).
    \]

    \item \textbf{Identity.}
    There is a distinguished element
    $\mathbf 1\in\mathcal O(1)$.

    \item \textbf{Symmetric group actions.}
    Each $\mathcal O(n)$ is equipped with a
    right action of the symmetric
    group $\Sigma_n$.
\end{enumerate}

These structures satisfy the following axioms:

\begin{enumerate}
    \item \textbf{Associativity.}
    Let $c\in\mathcal O(k)$,
    $d_s\in\mathcal O(n_s)$ for $1\leq s\leq k$,
    and $e_t\in\mathcal O(m_t)$ for
    $1\leq t\leq n$, where
    $n=n_1+\cdots+n_k$. Then
    \[
    \gamma\bigl(
        \gamma(c;d_1,\ldots,d_k);
        e_1,\ldots,e_n
    \bigr)
    =
    \gamma(c;f_1,\ldots,f_k),
    \]
    where
    \[
    f_s=
    \gamma\bigl(
        d_s;
        e_{N_{s-1}+1},\ldots,e_{N_s}
    \bigr),
    \qquad
    N_s=\sum_{r=1}^s n_r,\quad N_0=0.
    \]
    If $n_s=0$, the corresponding list is
    empty and $f_s=d_s$.

    \item \textbf{Unitality.}
    For all $c\in\mathcal O(n)$,
    \[
    \gamma(\mathbf 1;c)=c,
    \qquad
    \gamma(c;\mathbf 1,\ldots,\mathbf 1)=c.
    \]

    \item \textbf{Equivariance.}
    For $c\in\mathcal O(k)$,
    $d_s\in\mathcal O(n_s)$,
    $\sigma\in\Sigma_k$, and
    $\tau_s\in\Sigma_{n_s}$, we have
    \[
    \begin{aligned}
    &\gamma(c\sigma;d_1,\ldots,d_k)\\
    &\quad=
    \gamma(c;
        d_{\sigma^{-1}(1)},\ldots,
        d_{\sigma^{-1}(k)})
        \sigma(n_1,\ldots,n_k),
    \end{aligned}
    \]
    and
    \[
    \begin{aligned}
    &\gamma(c;d_1\tau_1,\ldots,d_k\tau_k)\\
    &\quad=
    \gamma(c;d_1,\ldots,d_k)
    (\tau_1\oplus\cdots\oplus\tau_k).
    \end{aligned}
    \]
    Here $\sigma(n_1,\ldots,n_k)$ denotes
    the corresponding block permutation,
    and $\tau_1\oplus\cdots\oplus\tau_k$
    denotes the permutation acting
    independently on each block.
\end{enumerate}
\end{definition}

We now explain how the structure of an operad naturally reflects the organization of abstract syntax trees in programming languages.

The elements of $\mathcal O(n)$ represent $n$-ary operations, which can be viewed as nodes with $n$ input slots in an abstract syntax tree. The composition map $\gamma$ describes how operations can be substituted into the input slots of other operations, corresponding to the recursive construction of complex expressions from simpler ones. Associativity ensures that, for a fixed tree, the resulting operation is independent of the order in which these substitutions are carried out. The identity element $\mathbf 1$ represents the operation that leaves its input unchanged.

The symmetric group actions describe the rearrangement of suboperations in an abstract syntax tree, while the equivariance axioms specify how the input arguments must be permuted accordingly to preserve the resulting composite operation. In particular, when the input slots of a node are reordered, each subtree moves together with its entire block of inputs. This induces the corresponding block permutation on the inputs of the composite.

The morphism between operads is a family of maps that preserves the operadic structure, including composition, identity, and symmetric group actions. The definition of a morphism between operads is as follows.

\begin{definition}
Let $\mathcal O$ and $\mathcal P$ be symmetric operads.
An \emph{operad morphism} $f:\mathcal O\to\mathcal P$
is a family of maps
\[
    f_n:\mathcal O(n)\longrightarrow\mathcal P(n),
    \qquad n\geq 0,
\]
satisfying the following conditions:

\begin{enumerate}
    \item \textbf{Compatibility with composition.}
    For all $k\geq 0$ and $n_1,\ldots,n_k\geq 0$,
    the following diagram commutes:
    \[
    \begin{tikzcd}
    \mathcal O(k)\times
    \mathcal O(n_1)\times\cdots\times\mathcal O(n_k)
    \arrow[r,"\gamma_{\mathcal O}"]
    \arrow[d,"f_k\times f_{n_1}\times\cdots\times f_{n_k}"']
    & \mathcal O(n_1+\cdots+n_k)
    \arrow[d,"f_{n_1+\cdots+n_k}"] \\
    \mathcal P(k)\times
    \mathcal P(n_1)\times\cdots\times\mathcal P(n_k)
    \arrow[r,"\gamma_{\mathcal P}"']
    & \mathcal P(n_1+\cdots+n_k).
    \end{tikzcd}
    \]

    \item \textbf{Preservation of the unit.}
    \[
        f_1(\mathbf 1_{\mathcal O})
        =\mathbf 1_{\mathcal P}.
    \]

    \item \textbf{Equivariance.}
    For every $n\geq 0$, $\theta\in\mathcal O(n)$,
    and $\sigma\in\Sigma_n$, we have
    \[
        f_n(\theta\cdot\sigma)
        =f_n(\theta)\cdot\sigma.
    \]
\end{enumerate}
\end{definition}
In general, a programming language may have multiple types, and the composition of operations must respect the types of their inputs and outputs. This leads to the notion of colored operads, where each operation has specified input and output colors corresponding to their types. The formal definition of colored operads is as follows.

\begin{definition}
A \emph{colored operad} $\mathcal O$ consists of
a set $C$ of \emph{colors} and a family of sets
\[
    \left\{
    \mathcal O(c_1,\ldots,c_n;c)
    \right\}_{n\geq 0,\;c_1,\ldots,c_n,c\in C},
\]
together with the following data:

\begin{enumerate}
    \item \textbf{Composition.}
    For every $k\geq 0$, colors
    $c,c_1,\ldots,c_k\in C$, and
    $d_{ij}\in C$ for
    $1\leq i\leq k$ and $1\leq j\leq n_i$,
    there is a map
    \[
    \begin{aligned}
    \gamma:\;
    &\mathcal O(c_1,\ldots,c_k;c)
    \times
    \prod_{i=1}^k
    \mathcal O(d_{i1},\ldots,d_{in_i};c_i)
    \\
    &\longrightarrow
    \mathcal O(
        d_{11},\ldots,d_{1n_1},
        \ldots,
        d_{k1},\ldots,d_{kn_k};c).
    \end{aligned}
    \]

    \item \textbf{Identity.}
    For each color $c\in C$, there is
    a distinguished element
    \[
        \mathbf 1_c\in\mathcal O(c;c).
    \]

    \item \textbf{Symmetric group actions.}
    For each $n\geq 0$ and
    $\sigma\in\Sigma_n$, there is a bijection
    \[
    \begin{aligned}
    \mathcal O(c_1,\ldots,c_n;c)
    &\longrightarrow
    \mathcal O(c_{\sigma(1)},\ldots,
               c_{\sigma(n)};c),\\
    f&\longmapsto f\sigma,
    \end{aligned}
    \]
    satisfying
    \[
        f\,\mathrm{id}=f,
        \qquad
        (f\sigma)\tau=f(\sigma\tau).
    \]
\end{enumerate}

These structures satisfy the following axioms:

\begin{enumerate}
    \item \textbf{Associativity.}
    Let
    \[
        f\in\mathcal O(c_1,\ldots,c_k;c),
    \]
    \[
        g_i\in
        \mathcal O(d_{i1},\ldots,d_{in_i};c_i),
    \]
    and
    \[
        h_{ij}\in
        \mathcal O(e_{ij1},\ldots,e_{ijm_{ij}};d_{ij}).
    \]
    Then
    \[
    \begin{aligned}
    &\gamma\bigl(
        \gamma(f;g_1,\ldots,g_k);
        h_{11},\ldots,h_{1n_1},
        \ldots,h_{k1},\ldots,h_{kn_k}
    \bigr)\\
    &\qquad=
    \gamma\bigl(
        f;
        \gamma(g_1;h_{11},\ldots,h_{1n_1}),
        \ldots,
        \gamma(g_k;h_{k1},\ldots,h_{kn_k})
    \bigr).
    \end{aligned}
    \]

    \item \textbf{Unitality.}
    For every
    $f\in\mathcal O(c_1,\ldots,c_n;c)$,
    \[
        \gamma(\mathbf 1_c;f)=f,
        \qquad
        \gamma(f;\mathbf 1_{c_1},\ldots,
                 \mathbf 1_{c_n})=f.
    \]

    \item \textbf{Equivariance.}
    Let
    \[
        f\in\mathcal O(c_1,\ldots,c_k;c),
    \]
    \[
        g_i\in
        \mathcal O(d_{i1},\ldots,d_{in_i};c_i).
    \]
    For $\sigma\in\Sigma_k$ and
    $\tau_i\in\Sigma_{n_i}$, we have
    \[
    \begin{aligned}
    &\gamma(
        f\sigma;
        g_{\sigma(1)},\ldots,g_{\sigma(k)}
    )\\
    &\qquad=
    \gamma(f;g_1,\ldots,g_k)
    \,\sigma\langle n_1,\ldots,n_k\rangle,
    \end{aligned}
    \]
    and
    \[
    \begin{aligned}
    &\gamma(
        f;
        g_1\tau_1,\ldots,g_k\tau_k
    )\\
    &\qquad=
    \gamma(f;g_1,\ldots,g_k)
    (\tau_1\oplus\cdots\oplus\tau_k).
    \end{aligned}
    \]
    Here
    $\sigma\langle n_1,\ldots,n_k\rangle$
    denotes the block permutation that
    rearranges consecutive blocks of lengths
    $n_1,\ldots,n_k$ into the order
    $\sigma(1),\ldots,\sigma(k)$, and
    $\tau_1\oplus\cdots\oplus\tau_k$
    acts independently on each block.
\end{enumerate}

Empty input lists and empty products are interpreted
in the usual way.
\end{definition}

A morphism of colored operads $\tau:\mathcal O\to\mathcal P$ consists of a map $\alpha$ between their sets of colors and, for each profile $(c_1,\ldots,c_n;c)$, a map
\[
    \mathcal O(c_1,\ldots,c_n;c)\longrightarrow
    \mathcal P(\alpha(c_1),\ldots,\alpha(c_n);\alpha(c)).
\]
These maps preserve composition, units, and the symmetric group actions.

\subsection{Algebras over Colored Operads}
An algebra over a colored operad $\mathcal O$ is a family of sets equipped with a family of operations corresponding to the elements of $\mathcal O$, together with a compatibility condition that ensures that the operations respect the composition and symmetry of the operad. The formal definition is as follows.
\begin{definition}
Let $\mathcal O$ be a symmetric $T$-colored operad in $\mathbf{Set}$.
An \emph{algebra over $\mathcal O$} is a family of sets
$A=(A_t)_{t\in T}$ equipped with structure maps
\[
    \theta_{t_1,\ldots,t_n;t}:
    \mathcal O(t_1,\ldots,t_n;t)
    \times A_{t_1}\times\cdots\times A_{t_n}
    \longrightarrow A_t,
    \qquad n\geq 0,
\]
for all $t_1,\ldots,t_n,t\in T$.
We write $\theta(o;\mathbf a)$ for the value of the
corresponding structure map. These maps satisfy the following conditions:

\begin{enumerate}
    \item Let $o\in\mathcal O(t_1,\ldots,t_k;t)$, and for each
    $1\leq i\leq k$, let
    $p_i\in\mathcal O(s_{i1},\ldots,s_{in_i};t_i)$ and
    $\mathbf a_i=(a_{i1},\ldots,a_{in_i})$ with
    $a_{ij}\in A_{s_{ij}}$.
    Write $n=n_1+\cdots+n_k$, and let $\mathbf a$ be the
    concatenation of $\mathbf a_1,\ldots,\mathbf a_k$.
    Then
    \[
        \theta\bigl(\gamma(o;p_1,\ldots,p_k);\mathbf a\bigr)
        =
        \theta\bigl(o;
            \theta(p_1;\mathbf a_1),\ldots,
            \theta(p_k;\mathbf a_k)\bigr).
    \]

    \item For every $t\in T$ and $a\in A_t$, the operadic unit
    $\mathbf 1_t\in\mathcal O(t;t)$ satisfies
    \[
        \theta(\mathbf 1_t;a)=a.
    \]

    \item For $o\in\mathcal O(t_1,\ldots,t_n;t)$,
    $\sigma\in\Sigma_n$, and $a_i\in A_{t_i}$, we have
    \[
        \theta\bigl(o\sigma;
            a_{\sigma(1)},\ldots,a_{\sigma(n)}\bigr)
        =\theta(o;a_1,\ldots,a_n).
    \]
\end{enumerate}
\end{definition}
\begin{example}
\label{example:free-construction}
Let $X$ be a program with a single type, and let $\mathcal O_X$ be the free symmetric operad generated by its chosen basic operation symbols.

First, introduce a formal symbol for each input and output under consideration, including the variables of $X$, and let $S_X$ be the resulting set of symbols. Form the free $\mathcal O_X$-algebra $F_X$: its elements are finite formal expressions built from these symbols and the basic operations, and its operations act by substitution of expressions.

Second, let $R_X$ be the original computational relations between these expressions. For instance, if $\widehat 0\in S_X$ names the value zero and $0+0=0$ is a given relation, include
\[
    +(\widehat 0,\widehat 0)=\widehat 0
\]
in $R_X$. Let $\sim_{R_X}$ be the smallest $\mathcal O_X$-algebra congruence containing $R_X$, and define
\[
    FA_X:=F_X/{\sim_{R_X}}.
\]
Thus, $FA_X$ consists of formal expressions modulo the original relations. If a constructed expression is identified with an original input or output symbol, its class is represented by that symbol; otherwise, its class is a new formal element. Its structure maps
\[
    \theta_n:\mathcal O_X(n)\times FA_X^n
    \longrightarrow FA_X
\]
are induced by substitution. If $f$ is binary and $g$ is unary, then the class represented by $f(g(x_1),x_2)$ satisfies
\[
    f(g(x_1),x_2)
    =\theta_2(f;\theta_1(g;x_1),x_2)
    =\theta_2(f\circ(g,1);x_1,x_2),
\]
where $1\in\mathcal O_X(1)$ is the identity operation and $f\circ(g,1)\in\mathcal O_X(2)$ is their operadic composite.
\end{example}

The construction has two steps: freely form expressions from the input and output symbols, and then impose their given relations by taking a quotient.

\begin{definition}[Morphism of Operadic Algebras]
\label{def:operadic-algebra-morphism}
Let $\mathcal O$ and $\mathcal P$ be colored operads with
color sets $T$ and $S$, respectively, and let
$A=(A_t)_{t\in T}$ be an $\mathcal O$-algebra and
$B=(B_s)_{s\in S}$ a $\mathcal P$-algebra.
A \emph{morphism of operadic algebras}
$(\tau,h):(\mathcal O,A)\to(\mathcal P,B)$ consists of
a morphism of colored operads
\[
    \tau:\mathcal O\longrightarrow\mathcal P
\]
with underlying color map $\alpha:T\to S$, together with
a family of maps
\[
    h_t:A_t\longrightarrow B_{\alpha(t)},\qquad t\in T,
\]
such that, for every $n\geq 0$ and $t_1,\ldots,t_n,t\in T$,
the following diagram commutes:
\[
\begin{tikzcd}[column sep=large, row sep=large]
\mathcal O(t_1,\ldots,t_n;t)\times\prod_{i=1}^n A_{t_i}
\arrow[r,"\theta^A"]
\arrow[d,"\tau\times\prod_{i=1}^n h_{t_i}"']
& A_t \arrow[d,"h_t"] \\
\mathcal P(\alpha(t_1),\ldots,\alpha(t_n);\alpha(t))
\times\prod_{i=1}^n B_{\alpha(t_i)}
\arrow[r,"\theta^B"']
& B_{\alpha(t)}.
\end{tikzcd}
\]
Here $\theta^A$ and $\theta^B$ are the corresponding algebra
structure maps, and $\tau$ on the left denotes its component
at $(t_1,\ldots,t_n;t)$.
When $n=0$, the empty products are the singleton set, so
the left vertical arrow is identified with the corresponding
nullary component of $\tau$.
For algebras with distinguished basepoints, we additionally
require $h_t(\ast)=\ast$ for every $t\in T$.
\end{definition}

The semantics of programming languages can be understood in terms of the input-output behavior of programs: given an input to a program, what output does it produce? Operadic algebras provide a natural framework for modeling this behavior, as they capture the compositional structure of operations and their interactions with inputs and outputs. In particular, the homomorphism of operadic algebras allows us to relate the semantics of different programming languages by mapping operations and their corresponding inputs and outputs from one language to another. The square in Definition~\ref{def:operadic-algebra-morphism} says that evaluating a source operation and then mapping its result gives the same value as mapping its inputs and evaluating the translated operation.

Suppose that we have a programming language $\mathcal{L}$ with a set of types $T$. For each finite sequence of types $(t_1,\ldots,t_n)$ and each $t\in T$, let $\Sigma(t_1,\ldots,t_n;t)$ denote the set of basic operation symbols that take inputs of types $t_1,\ldots,t_n$ and produce an output of type $t$. These sets form a $T$-colored signature $\Sigma$. We write $\mathcal O_{\mathcal L}$ for the free colored operad generated by this chosen signature. Its colors are the types in $T$, and its operations are obtained by composing the basic operation symbols whenever their input and output types match.

An operation on ordinary values may be undefined for some arguments. Write $f_A$ for the interpretation of a basic operation $f$ on ordinary values. If operations with interpretations $g_1,\ldots,g_k$ are substituted into the input slots of $f$, we assume that the resulting expression in $\mathcal L$ is evaluated by
\begin{equation}
\label{eq:language-substitution}
    \bigl(f_A\circ(g_1,\ldots,g_k)\bigr)(\mathbf a)
    = f_A\bigl(g_1(\mathbf a_1),\ldots,g_k(\mathbf a_k)\bigr),
\end{equation}
where $\mathbf a$ is the concatenation of the input blocks $\mathbf a_1,\ldots,\mathbf a_k$. The two sides are defined on exactly the same input tuples.

We consider two related constructions. The first freely forms expressions from input and output symbols and retains their given relations (see Example~\ref{example:free-construction}). The second works directly with value domains and uses a basepoint to represent undefined outcomes. We begin with the formal construction.

First, for every type $t\in T$, let $S_t$ be a set of formal symbols naming all inputs and outputs under consideration. Write $S=(S_t)_{t\in T}$ and form the free $\mathcal O_{\mathcal L}$-algebra
\[
    F_{\mathcal L}=(F_{\mathcal L,t})_{t\in T}.
\]
Its elements are finite well-typed terms: a symbol $s\in S_t$ is a term of type $t$, and if $f\in\Sigma(t_1,\ldots,t_n;t)$ and $u_i$ is a term of type $t_i$, then $f(u_1,\ldots,u_n)$ is a term of type $t$. This includes nullary operations. An operation $o\in\mathcal O_{\mathcal L}(t_1,\ldots,t_n;t)$ acts by substituting the terms $u_i$ into the inputs of its expression tree, giving structure maps
\[
    \theta^{F_{\mathcal L}}:
    \mathcal O_{\mathcal L}(t_1,\ldots,t_n;t)
    \times F_{\mathcal L,t_1}\times\cdots\times F_{\mathcal L,t_n}
    \longrightarrow F_{\mathcal L,t}.
\]

Second, let $R$ be the prescribed relations between terms with the same output type. In particular, a given operation with input symbols $s_1,\ldots,s_n$ and output symbol $s$ contributes the relation
\[
    f(s_1,\ldots,s_n)=s.
\]
Let $\sim_R=(\sim_{R,t})_{t\in T}$ be the smallest $\mathcal O_{\mathcal L}$-algebra congruence containing $R$. Define
\[
    FA_{\mathcal L}:=F_{\mathcal L}/{\sim_R}
    =\left(F_{\mathcal L,t}/{\sim_{R,t}}\right)_{t\in T}.
\]
Write $FA_{\mathcal L,t}$ for its component of type $t$ and $q_t:F_{\mathcal L,t}\to FA_{\mathcal L,t}$ for the quotient map. The induced structure maps are defined by
\[
    \theta^{FA_{\mathcal L}}(o;q_{t_1}(u_1),\ldots,q_{t_n}(u_n))
    :=q_t\bigl(\theta^{F_{\mathcal L}}(o;u_1,\ldots,u_n)\bigr).
\]

\begin{proposition}
\label{prop:free-operadic-algebra}
The family $FA_{\mathcal L}$, equipped with the induced structure maps $\theta^{FA_{\mathcal L}}$, is an $\mathcal O_{\mathcal L}$-algebra.
\end{proposition}

\begin{proof}
If $q_{t_i}(u_i)=q_{t_i}(v_i)$ for every $i$, then $u_i\sim_R v_i$. Since $\sim_R$ is an algebra congruence, $\theta^{F_{\mathcal L}}(o;u_1,\ldots,u_n)\sim_R\theta^{F_{\mathcal L}}(o;v_1,\ldots,v_n)$. Hence the induced structure maps $\theta^{FA_{\mathcal L}}$ are well-defined on the quotient algebra $FA_{\mathcal L}$.

Let $(a_1,\ldots,a_n)$ be a tuple of inputs in the quotient algebra $FA_{\mathcal L}$. For each $j$, choose a representative $u_j\in F_{\mathcal L,t_j}$ of $a_j$. Writing $\mathbf u=(u_1,\ldots,u_n)$, we then have $q\mathbf u=(a_1,\ldots,a_n)$. Let $o,p_1,\ldots,p_k$ be composable operations whose composite has input types $t_1,\ldots,t_n$. Partition $\mathbf u$ into the consecutive input blocks $\mathbf u_1,\ldots,\mathbf u_k$ for $p_1,\ldots,p_k$. By the definition of the quotient structure maps and the substitution law in $F_{\mathcal L}$, we have
\[
\begin{aligned}
    \theta^{FA_{\mathcal L}}(\gamma(o;p_1,\ldots,p_k);q\mathbf u)
    &=q\bigl(\theta^{F_{\mathcal L}}(\gamma(o;p_1,\ldots,p_k);\mathbf u)\bigr)\\
    &=q\bigl(\theta^{F_{\mathcal L}}(o;
        \theta^{F_{\mathcal L}}(p_1;\mathbf u_1),\ldots,
        \theta^{F_{\mathcal L}}(p_k;\mathbf u_k))\bigr)\\
    &=\theta^{FA_{\mathcal L}}(o;\theta^{FA_{\mathcal L}}(p_1;q\mathbf u_1),\ldots,
        \theta^{FA_{\mathcal L}}(p_k;q\mathbf u_k)).
\end{aligned}
\]
The unit law follows from
\[
    \theta^{FA_{\mathcal L}}(\mathbf 1_t;q_t(u))
    =q_t(\theta^{F_{\mathcal L}}(\mathbf 1_t;u))=q_t(u).
\]
For an $n$-ary operation $o$ and $\sigma\in\Sigma_n$, write $\mathbf u^\sigma=(u_{\sigma(1)},\ldots,u_{\sigma(n)})$. Equivariance follows from
\[
\begin{aligned}
    \theta^{FA_{\mathcal L}}(o\sigma;(q\mathbf u)^\sigma)
    &=q\bigl(\theta^{F_{\mathcal L}}(o\sigma;\mathbf u^\sigma)\bigr)\\
    &=q\bigl(\theta^{F_{\mathcal L}}(o;\mathbf u)\bigr)
      =\theta^{FA_{\mathcal L}}(o;q\mathbf u).
\end{aligned}
\]
Every element of $FA_{\mathcal L,t}$ has a representative under $q_t$, so these identities verify all the operadic algebra axioms on the quotient algebra $FA_{\mathcal L}$.
\end{proof}

The free algebra $F_{\mathcal L}$ contains every finite well-typed expression that can be built from the symbols in $S$ using the operations in $\Sigma$. It may therefore contain formal composites beyond those explicitly given in $\mathcal L$. To ensure that these expressions respect the prescribed relations $R$ of the original language, we take the quotient of $F_{\mathcal L}$ by the congruence generated by $R$. In the resulting quotient algebra $FA_{\mathcal L}$, a relation such as $f(a)=b$ allows us to replace $f(a)$ by $b$ in any well-typed expression without changing its class.

\begin{definition}[Free Semantic Algebra]
\label{def:free-semantic-algebra}
The $\mathcal O_{\mathcal L}$-algebra $FA_{\mathcal L}$ constructed above is called the \emph{free semantic algebra of the programming language $\mathcal L$ on the input and output symbols $S$ subject to the given relations $R$}.
\end{definition}

An assignment of the generators $S$ to elements of the corresponding types in another $\mathcal O_{\mathcal L}$-algebra extends uniquely to an algebra morphism from $FA_{\mathcal L}$ precisely when it respects the relations $R$.

We now give the second construction using the basic operations on ordinary values.

For each $t\in T$, let $A_t$ be the ordinary value domain of type $t$. Choose a fresh element $\ast$, called the \emph{undefined outcome}, and let $\mathsf{Undef}=\{\ast\}$ be its distinguished outcome type. Set
\[
    A_t^\ast=A_t\cup\{\ast\}.
\]
For each basic operation $f\in\Sigma(t_1,\ldots,t_n;t)$, its interpretation on ordinary values has the form
\[
    f_A:\operatorname{dom}(f_A)\longrightarrow A_t,\qquad
    \operatorname{dom}(f_A)\subseteq A_{t_1}\times\cdots\times A_{t_n}.
\]
Its strict totalization is the map
\[
    f^\ast:A_{t_1}^\ast\times\cdots\times A_{t_n}^\ast
    \longrightarrow A_t^\ast
\]
defined by
\[
    f^\ast(a_1,\ldots,a_n)=
    \begin{cases}
        f_A(a_1,\ldots,a_n),
        &\text{if all $a_i\in A_{t_i}$ and $f_A(a_1,\ldots,a_n)$ is defined},\\
        \ast,&\text{otherwise}.
    \end{cases}
\]
Thus, undefined operations produce $\ast$, and every operation propagates an undefined argument to its result.

To evaluate an operation of $\mathcal O_{\mathcal L}$ on $A^\ast=(A_t^\ast)_{t\in T}$, place the inputs on the corresponding edges of its tree and apply $f^\ast$ at each vertex labeled by $f$. The value on the output edge defines the structure maps
\[
    \theta^{A^\ast}:
    \mathcal O_{\mathcal L}(t_1,\ldots,t_n;t)
    \times A_{t_1}^\ast\times\cdots\times A_{t_n}^\ast
    \longrightarrow A_t^\ast.
\]

\begin{proposition}
\label{prop:semantic-algebra}
The family $A^\ast$, equipped with the structure maps $\theta^{A^\ast}$, is an $\mathcal O_{\mathcal L}$-algebra.
\end{proposition}

\begin{proof}
Let $o,p_1,\ldots,p_k$ be composable operations and let $\mathbf a$ be the concatenation of the input tuples $\mathbf a_i$ for $p_i$. We verify the substitution law
\[
    \theta^{A^\ast}\bigl(\gamma(o;p_1,\ldots,p_k);\mathbf a\bigr)
    =\theta^{A^\ast}\bigl(o;
        \theta^{A^\ast}(p_1;\mathbf a_1),\ldots,
        \theta^{A^\ast}(p_k;\mathbf a_k)\bigr).
\]
Since $\mathcal O_{\mathcal L}$ is freely generated by $\Sigma$, it suffices to verify the substitution law when the outer operation is basic. Take $o=f\in\Sigma$, and let $g_i$ be the interpretation of $p_i$ on ordinary inputs. We claim that
\begin{equation}
\label{eq:strict-totalization-substitution}
    \bigl(f_A\circ(g_1,\ldots,g_k)\bigr)^\ast(\mathbf a)
    = f^\ast\bigl(g_1^\ast(\mathbf a_1),\ldots,g_k^\ast(\mathbf a_k)\bigr).
\end{equation}
By~\eqref{eq:language-substitution}, if the composite is defined at $\mathbf a$, each $g_i$ is defined on its input block and $f_A$ is defined on the resulting tuple. Both sides of~\eqref{eq:strict-totalization-substitution} then equal $f_A(g_1(\mathbf a_1),\ldots,g_k(\mathbf a_k))$.

Suppose now that the left-hand side of~\eqref{eq:strict-totalization-substitution} is $\ast$. If some $g_i^\ast(\mathbf a_i)=\ast$, strictness of $f^\ast$ gives $\ast$ on the right as well. Otherwise each $b_i=g_i^\ast(\mathbf a_i)$ is ordinary and equals $g_i(\mathbf a_i)$. If $f_A$ were defined at $(b_1,\ldots,b_k)$, then~\eqref{eq:language-substitution} would make the composite defined at $\mathbf a$, contrary to our assumption. Hence $f^\ast(b_1,\ldots,b_k)=\ast$, proving~\eqref{eq:strict-totalization-substitution}.

Starting at the input edges and applying~\eqref{eq:strict-totalization-substitution} at each vertex, we obtain the strict totalization of the operation represented by the whole tree. In the tree for $\gamma(o;p_1,\ldots,p_k)$ we may therefore evaluate the subtrees $p_i$ first and pass their values to $o$, which gives the substitution law above. The identity on $A_t$ extends to the identity on $A_t^\ast$ because its extension also fixes $\ast$. Finally, relabeling the input slots and their arguments together preserves both definedness and the resulting value, so totalization also commutes with input permutations.
\end{proof}

\begin{definition}[Semantic Algebra]
\label{def:semantic-algebra}
The $\mathcal O_{\mathcal L}$-algebra $A^\ast$ constructed above is called a \emph{semantic algebra of the programming language $\mathcal L$}.
\end{definition}

Even when every basic operation is defined on all tuples
of ordinary values and returns an ordinary value, we may
adjoin a fresh basepoint $\ast$ to each value domain and
extend the operations strictly as above.
This extension preserves the original operations on
ordinary values and returns $\ast$ whenever an input is
$\ast$. Thus, it introduces no undefined outcomes on
ordinary inputs and still yields an $\mathcal O_{\mathcal L}$-algebra.
The same construction applies to free semantic algebras.
Henceforth, all semantic algebras are understood to carry
a distinguished basepoint in each type, whether or not
their interpretations on ordinary values exhibit
undefined behavior.

\begin{remark}
\label{rem:program-semantic-algebras}
To specialize Definitions~\ref{def:free-semantic-algebra} and~\ref{def:semantic-algebra} to a program $X$, choose its typed operation signature $\Sigma_X$ and corresponding operad $\mathcal O_X$. For the free construction, use its input and output symbols $S_X$ and their original relations $R_X$. For the value-based construction, use its value domains and the interpretations of its basic operations on those domains. We denote the resulting free semantic algebra and semantic algebra by $FA_X$ and $A_X$, respectively. Unless otherwise stated, these names refer to the two constructions above, with the preceding basepoint convention applied when needed.
\end{remark}

The semantic algebra constructed above captures how
undefined behavior can make the order of addition
semantically significant. For mathematical integers,
associativity allows the additions in a sum to be performed
in different orders without changing its value.
In a programming language, however, one order of addition
may encounter an undefined
intermediate operation while another remains defined.
We therefore cannot assume that changing the order of
addition preserves its semantics.

For example, let $M$ be the largest value of a signed
integer type. If overflowing addition is interpreted as
$\ast$, then
\[
    (M+1)+(-1)=\ast,
    \qquad
    M+(1+(-1))=M.
\]
The first expression contains the undefined intermediate
addition $M+1$, and $\ast$ propagates to the final result.
The second first computes $1+(-1)=0$ and then returns $M$.
Thus, although both orders of addition yield $M$ over the
mathematical integers, they have different meanings in
this semantic algebra. Consequently, the ordinary
associativity equation for addition cannot be imposed
in this model.

This distinction is important for our \Cxx{}-to-\Lean{} compiler. The semantics of \Cxx{} include undefined behavior, so changing the order of addition can change whether a computation is defined, as the example above shows. Our compiler must therefore express the source order of addition explicitly in the generated \Lean{} code. Mathematical associativity alone does not justify changing that order; such a transformation must account for undefined behavior and preserve the intended source semantics.

We now compare the free semantic algebra $FA_{\mathcal L}$ of Definition~\ref{def:free-semantic-algebra} with the semantic algebra $A_{\mathcal L}=(A_t^\ast)_{t\in T}$ of Definition~\ref{def:semantic-algebra}, using the same chosen input and output data. Here $S_t$ consists of the symbols $\widehat a$, one for each semantic value $a\in A_t^\ast$, together with a disjoint family of variable symbols $V_t$. The relations $R$ consist of all concrete computation relations on these values, including computations with result $\ast$, and the basepoint propagation relations, whose other arguments may be arbitrary well-typed terms. The class of $\widehat\ast_t$ is the distinguished basepoint. We use the free semantic algebra already constructed above.

\begin{proposition}
\label{prop:semantic-restriction-and-embedding}
There are basepoint-preserving operadic algebra morphisms over $\mathrm{id}_{\mathcal O_{\mathcal L}}$
\[
    \operatorname{res}:FA_{\mathcal L}\twoheadrightarrow A_{\mathcal L},
    \qquad
    \iota:A_{\mathcal L}\hookrightarrow FA_{\mathcal L},
\]
satisfying $\operatorname{res}\circ\iota=\mathrm{id}_{A_{\mathcal L}}$.
\end{proposition}

\begin{proof}
Let $z\in(FA_{\mathcal L})_t$ and choose a representative $u$. If $u$ has no defined value in $\mathcal L$ on the given data, set $\operatorname{res}_t(z)=\ast$. If $u$ has value $a\in A_t$, set $\operatorname{res}_t(z)=a$. Here each formal variable is assigned $\ast$.

If $v$ is another representative of $z$, then $u\sim_R v$. Both sides of each concrete computation relation in $R$ have the same outcome in $A_{\mathcal L}$, and both sides of each propagation relation give $\ast$ by strictness. Replacing one side of such a relation by the other therefore preserves the outcome in any surrounding expression. Thus $u$ and $v$ have the same outcome in $A_{\mathcal L}$, so $\operatorname{res}_t$ is well-defined.

For $a\in(A_{\mathcal L})_t$, let $\iota_t(a)$ be the class of the generator corresponding to $a$. The concrete computation relations ensure that $\iota$ preserves each basic operation. To check compatibility for $\operatorname{res}$, choose representatives $u_i$ of the input classes $z_i$. Then $f(u_1,\ldots,u_n)$ represents the result of applying $f$ to those classes, and its value is $f^\ast(\operatorname{res}_{t_1}(z_1),\ldots,\operatorname{res}_{t_n}(z_n))$. This also covers nullary operations.

By the algebra axioms, the operations preserved by both maps include the units and are closed under substitution and input permutations. Since they include the generators $\Sigma$, both maps preserve every operation of $\mathcal O_{\mathcal L}$.

By construction, $\iota_t(\ast)$ is the distinguished basepoint and $\operatorname{res}_t(\iota_t(a))=a$ for every $a$, including $\ast$. Hence both maps preserve the basepoint, $\operatorname{res}\circ\iota=\mathrm{id}_{A_{\mathcal L}}$, and $\operatorname{res}$ is surjective and $\iota$ is injective in each type.
\end{proof}

\begin{remark}[Geometric intuition]
\label{rem:semantic-deformation}
From a geometric viewpoint, $FA_{\mathcal L}$ may be regarded as a \emph{universal deformation} of $A_{\mathcal L}$. This is analogous in spirit to the relationship between a space and its tangent spaces: just as tangent spaces describe possible directions of variation, the operad specifies the directions in which the program can be extended.
\end{remark}

Inspired by Proposition~\ref{prop:semantic-restriction-and-embedding}, we give the following definition.

\begin{definition}[Semantic Faithfulness and Semantic Restriction]
\label{def:semantic-faithfulness-and-restriction}
Let $(\tau,h):(\mathcal O,A)\to(\mathcal P,B)$ be a morphism of semantic algebras as in Definition~\ref{def:operadic-algebra-morphism}, with underlying type map $\alpha:T\to S$. We call it \emph{semantically faithful} if $h_t:A_t\to B_{\alpha(t)}$ is injective for every $t\in T$, and a \emph{semantic restriction} if $h_t$ is surjective for every $t\in T$.
\end{definition}

\begin{remark}
\label{rem:faithfulness-preserves-definedness}
A semantically faithful morphism $(\tau,h)$ preserves definedness: if $a\in A_t$ and $a\neq\ast$, then $h_t(a)\neq\ast$.
\end{remark}

For a program $X$ written in $\mathcal L$, the inclusions of its types and basic operations induce an operad embedding $j:\mathcal O_X\hookrightarrow\mathcal O_{\mathcal L}$ with injective type map $\alpha:T_X\hookrightarrow T$. We specify the interpretation of program inputs and outputs as language values by basepoint-preserving injections
\[
    e_t:(A_X)_t\hookrightarrow(A_{\mathcal L})_{\alpha(t)},
    \qquad t\in T_X.
\]
Use the presentations of Proposition~\ref{prop:semantic-restriction-and-embedding} for $FA_X$ and $FA_{\mathcal L}$, and choose generator assignments
\[
    \lambda_t:S_{X,t}\longrightarrow(FA_{\mathcal L})_{\alpha(t)},
    \qquad
    \lambda_t(\widehat a)=\widehat{e_t(a)}
    \quad(a\in(A_X)_t).
\]
Write $\lambda^\#$ for the extension to $F_X$ obtained by translating the basic operation symbols through $j$. Require every defining relation $u=v$ in $R_X$ of type $t$ to satisfy
\[
    \lambda^\#_t(u)=\lambda^\#_t(v)
    \quad\text{in }(FA_{\mathcal L})_{\alpha(t)}.
\]
Thus all relations of $FA_X$ are preserved in $FA_{\mathcal L}$. In particular, for every basic operation $f\in\Sigma_X(t_1,\ldots,t_n;t)$ and correctly typed inputs $a_i$, the computation relations give
\[
    e_t\bigl(\theta^{A_X}(f;a_1,\ldots,a_n)\bigr)
    =\theta^{A_{\mathcal L}}
      (j(f);e_{t_1}(a_1),\ldots,e_{t_n}(a_n)).
\]

The \emph{extension by zero} of $A_X$ along $j$, denoted by $j_!A_X$, has value sets
\[
    (j_!A_X)_s=
    \begin{cases}
        e_t((A_X)_t),&\text{if $s=\alpha(t)$ for some $t\in T_X$},\\
        \{\ast\},&\text{if $s\notin\alpha(T_X)$}.
    \end{cases}
\]
Require $j_!A_X$ to be closed under every basic language operation $f\in\Sigma(s_1,\ldots,s_n;s)$, in the sense that
\[
    \theta^{A_{\mathcal L}}(f;b_1,\ldots,b_n)\in(j_!A_X)_s
    \qquad\text{whenever }b_i\in(j_!A_X)_{s_i}.
\]
This requirement includes nullary operations. The extension retains the language computation relations by restricting the structure maps of $A_{\mathcal L}$ to $j_!A_X$,
\[
    \theta^{j_!A_X}(o;\mathbf b)
    :=\theta^{A_{\mathcal L}}(o;\mathbf b),
    \qquad b_i\in(j_!A_X)_{s_i},
\]
for $o\in\mathcal O_{\mathcal L}(s_1,\ldots,s_n;s)$. Write $\iota:j_!A_X\hookrightarrow A_{\mathcal L}$ for the inclusion.

Define basepoint-preserving surjections
\[
    q_s:(A_{\mathcal L})_s\longrightarrow(j_!A_X)_s,
    \qquad
    q_s(b)=
    \begin{cases}
        b,&\text{if $b\in(j_!A_X)_s$},\\
        \ast,&\text{otherwise}
    \end{cases}.
\]
Require the language computation relations to remain valid after replacing every value by its image under $q$: for each basic operation $f\in\Sigma(s_1,\ldots,s_n;s)$ and correctly typed values $b_1,\ldots,b_n,d$,
\begin{equation}
\label{eq:program-restriction-relations}
    \theta^{A_{\mathcal L}}(f;b_1,\ldots,b_n)=d
    \quad\Longrightarrow\quad
    \theta^{A_{\mathcal L}}
      (f;q_{s_1}(b_1),\ldots,q_{s_n}(b_n))=q_s(d).
\end{equation}
This requirement includes nullary operations.

\begin{proposition}
\label{prop:program-semantic-restriction}
With the restricted structure maps, $j_!A_X$ is an $\mathcal O_{\mathcal L}$-algebra. The inclusion $\iota:j_!A_X\hookrightarrow A_{\mathcal L}$ is semantically faithful. The maps $e_t$ give an isomorphism of $\mathcal O_X$-algebras from $A_X$ to the restriction of $j_!A_X$ along $j$. Moreover, $q\circ\iota=\mathrm{id}_{j_!A_X}$, and $q$ is a semantic restriction.
\end{proposition}

\begin{proof}
Closure under the basic language operations gives closure under every operadic composite. The structure maps therefore restrict from $A_{\mathcal L}$ and satisfy composition, units, and equivariance. Thus $j_!A_X$ is an $\mathcal O_{\mathcal L}$-algebra.

The inclusion $\iota$ preserves the structure maps by their definition as restrictions of those on $A_{\mathcal L}$. It preserves the basepoint and is injective in every type, so it is semantically faithful.

For each basic program operation $f$, the compatibility established above and the definition of the structure maps on $j_!A_X$ give
\[
    e_t\bigl(\theta^{A_X}(f;\mathbf a)\bigr)
    =\theta^{j_!A_X}(j(f);e\mathbf a).
\]
Here $e\mathbf a=(e_{t_1}(a_1),\ldots,e_{t_n}(a_n))$. This identity holds for the units and is preserved by operadic substitution and input permutations. Since the basic program operations generate $\mathcal O_X$, it holds for every operation in $\mathcal O_X$. At each program type $t$, the restriction along $j$ has value set $(j_!A_X)_{\alpha(t)}=e_t((A_X)_t)$. Each $e_t$ is a bijection onto this set, so $e$ gives the asserted isomorphism.

Since $q$ fixes every element of $j_!A_X$, we have $q\circ\iota=\mathrm{id}_{j_!A_X}$. For each basic language operation $f$, condition~\eqref{eq:program-restriction-relations} and the definition of the restricted structure maps give
\[
    q_s\bigl(\theta^{A_{\mathcal L}}(f;b_1,\ldots,b_n)\bigr)
    =\theta^{j_!A_X}
      (f;q_{s_1}(b_1),\ldots,q_{s_n}(b_n)).
\]
This compatibility extends to all operadic operations by composition, units, and symmetry. Hence $q$ is an operadic algebra morphism; its components are surjective and preserve the basepoints, so it is a semantic restriction.
\end{proof}

\begin{remark}
\label{rem:extension-by-zero}
Identify $A_X$ with its image under $e$. For every subalgebra $B\subseteq A_{\mathcal L}$ containing the basepoint of each type,
\[
    j_!A_X\subseteq B
    \quad\Longleftrightarrow\quad
    A_X\subseteq j^*B,
\]
where $j^*B$ denotes restriction along $j$. The unit of the adjunction $j_!\dashv j^*$ at $A_X$ is the isomorphism $A_X\xrightarrow{\sim}j^*j_!A_X$ of Proposition~\ref{prop:program-semantic-restriction}.
\end{remark}

In fact, this construction also helps us understand when a fragment of a program can be represented by a function. Such a representation requires the types occurring in the function to match the corresponding types already present in the program, and every operation used by the function to be defined in the program.

Suppose that a fragment of a program $X$ can be represented by a function $f$. We require the types and basic operation symbols used by $f$ to embed into those of $X$. These embeddings induce an operad embedding $j:\mathcal O_f\hookrightarrow\mathcal O_X$ with injective type map $\alpha:T_f\hookrightarrow T_X$. We specify the interpretation of the inputs and outputs of $f$ in the program by basepoint-preserving injections
\[
    e_t:(A_f)_t\hookrightarrow(A_X)_{\alpha(t)},
    \qquad t\in T_f.
\]
The maps $e_t$ determine the concrete parameter substitution: an input value $a\in(A_f)_t$ is assigned the program value $e_t(a)$. To preserve the semantics of $f$, we require that, for every $o\in\mathcal O_f(t_1,\ldots,t_n;t)$ and correctly typed inputs $a_i$,
\[
    e_t\bigl(\theta^{A_f}(o;a_1,\ldots,a_n)\bigr)
    =\theta^{A_X}
      (j(o);e_{t_1}(a_1),\ldots,e_{t_n}(a_n)).
\]

In the next section, we apply the construction of Proposition~\ref{prop:program-semantic-restriction} with a source language in place of the program $X$ and a target language in place of the ambient language $\mathcal L$.

\section{Compiler}
\label{sec:compiler}
Compiler is a machine that translates a program written in one programming language into another programming language. The translation is done in such a way that the semantics of the original program are preserved in the translated program. In this section, we define a compiler as a morphism of semantic algebras, and we show that such a morphism preserves the semantics of the original program.

We have already seen that we can translate semantic algebra of a program $X$ which is an algebra over a colored operad $\mathcal O_X$ into a semantic algebra of a programming language $\mathcal L$ which is an algebra over a colored operad $\mathcal O_{\mathcal L}$. This consturction is a general construction that can be applied to transalte semantic algebra of a program language into semantic algebra of another programming language. 

Let $L_1$ and $L_2$ be programming languages with associated colored operads $\mathcal O_{L_1}$ and $\mathcal O_{L_2}$ and semantic algebras $A_{L_1}$ and $A_{L_2}$, respectively. We construct a compiler from $L_1$ to $L_2$ using the extension by zero of Remark~\ref{rem:extension-by-zero}.

Let $T_i$ be the type set and $\Sigma_i$ the chosen basic operation signature of $L_i$, for $i=1,2$. We require the types and basic operation symbols of $L_1$ to embed into those of $L_2$, inducing an operad embedding with injective type map
\[
    j:\mathcal O_{L_1}\hookrightarrow\mathcal O_{L_2},
    \qquad
    \alpha:T_1\hookrightarrow T_2.
\]
Specify the interpretation of source inputs and outputs as target values by basepoint-preserving injections
\[
    e_t:(A_{L_1})_t\hookrightarrow(A_{L_2})_{\alpha(t)},
    \qquad t\in T_1.
\]
Require the chosen correspondences of generators and basic operations to send every defining relation of $FA_{L_1}$ to a relation valid in $FA_{L_2}$, with operations translated by $j$ and each value symbol $\widehat a$ sent to $\widehat{e_t(a)}$. In particular, for every basic operation $f\in\Sigma_1(t_1,\ldots,t_n;t)$ and correctly typed inputs $a_i$,
\[
    e_t\bigl(\theta^{A_{L_1}}(f;a_1,\ldots,a_n)\bigr)
    =\theta^{A_{L_2}}
      (j(f);e_{t_1}(a_1),\ldots,e_{t_n}(a_n)).
\]
This compatibility extends to all source operations, so $(j,e)$ is the forward morphism of semantic algebras.

Using $e$ to identify source values with their images in $A_{L_2}$, define the extension by zero $j_!A_{L_1}$ on each target type by
\[
    (j_!A_{L_1})_s=
    \begin{cases}
        e_t((A_{L_1})_t),&\text{if $s=\alpha(t)$ for some $t\in T_1$},\\
        \{\ast\},&\text{if $s\notin\alpha(T_1)$}.
    \end{cases}
\]
Require $j_!A_{L_1}$ to be closed under every basic target operation $f\in\Sigma_2(s_1,\ldots,s_n;s)$, in the sense that
\[
    \theta^{A_{L_2}}(f;b_1,\ldots,b_n)\in(j_!A_{L_1})_s
    \qquad\text{whenever }b_i\in(j_!A_{L_1})_{s_i}.
\]
This requirement includes nullary operations.

To describe the extension by generators and relations, form finite well-typed target-language expressions with entries from $j_!A_{L_1}$, and let $\operatorname{ev}_{2,s}$ denote their interpretation in $(A_{L_2})_s$. Closure ensures that every such interpretation lies in $(j_!A_{L_1})_s$. Impose the computation relations
\[
    u\sim_s^{\mathrm{ev}}v
    \quad\Longleftrightarrow\quad
    \operatorname{ev}_{2,s}(u)=\operatorname{ev}_{2,s}(v)
\]
for expressions of the same type. Their quotient is identified with $j_!A_{L_1}$, since each value already occurs as a generator. The structure maps are
\[
    \theta^{j_!A_{L_1}}(o;\mathbf b)
    :=\theta^{A_{L_2}}(o;\mathbf b),
    \qquad b_i\in(j_!A_{L_1})_{s_i},
\]
for $o\in\mathcal O_{L_2}(s_1,\ldots,s_n;s)$.

The inclusion $\iota:j_!A_{L_1}\hookrightarrow A_{L_2}$ preserves the restricted structure maps and is injective in every type, so it is semantically faithful. Restricting along $j$ recovers the source semantics. At type $t$ the value set is $e_t((A_{L_1})_t)$, and each operation $o\in\mathcal O_{L_1}$ acts through its image $j(o)\in\mathcal O_{L_2}$. With their codomains restricted to these images, the maps $e_t$ therefore give an isomorphism $A_{L_1}\xrightarrow{\sim}j^*j_!A_{L_1}$ by the computation compatibility above.

Define the reverse basepoint-preserving surjections
\[
    q_s:(A_{L_2})_s\longrightarrow(j_!A_{L_1})_s,
    \qquad
    q_s(b)=
    \begin{cases}
        b,&\text{if $b\in(j_!A_{L_1})_s$},\\
        \ast,&\text{otherwise}
    \end{cases}.
\]
These satisfy $q\circ\iota=\mathrm{id}_{j_!A_{L_1}}$. As in Proposition~\ref{prop:program-semantic-restriction}, they give a semantic restriction provided that every target-language computation relation remains valid after replacing its values by their images under $q$. Explicitly, for every basic operation $f\in\Sigma_2(s_1,\ldots,s_n;s)$ and correctly typed values $b_1,\ldots,b_n,d$, require
\[
    \theta^{A_{L_2}}(f;b_1,\ldots,b_n)=d
    \quad\Longrightarrow\quad
    \theta^{A_{L_2}}
      (f;q_{s_1}(b_1),\ldots,q_{s_n}(b_n))=q_s(d).
\]
This condition includes nullary operations and gives the required compatibility with every operadic operation.

\begin{definition}[Compiler]
\label{def:compiler}
The \emph{compiler from $L_1$ to $L_2$ determined by $(j,e)$} is the morphism of operadic algebras
\[
    (j,c):(\mathcal O_{L_1},A_{L_1})
    \longrightarrow
    (\mathcal O_{L_2},j_!A_{L_1})
\]
in which $j$ translates operations and $c:A_{L_1}\to j_!A_{L_1}$ translates values over the type map $\alpha$. Its component at each source type $t$ is
\[
    c_t:(A_{L_1})_t\longrightarrow(j_!A_{L_1})_{\alpha(t)},
    \qquad c_t(a)=e_t(a).
\]
Composing with the inclusion $\iota:j_!A_{L_1}\hookrightarrow A_{L_2}$ recovers the interpretation in the target language. The identities $\iota_{\alpha(t)}\circ c_t=e_t$ determine $c$ uniquely.
\end{definition}

We now construct the maps $e_t$ through the free semantic algebras. For $i=1,2$, use the presentations of $FA_{L_i}$ from Proposition~\ref{prop:semantic-restriction-and-embedding}, with generators $S_i$ containing a symbol $\widehat a$ for every semantic value, including the basepoint, and with the computation relations $R_i$ specified there. Choose a type correspondence $\alpha:T_1\hookrightarrow T_2$ and a correspondence of basic operation symbols inducing $j:\mathcal O_{L_1}\hookrightarrow\mathcal O_{L_2}$ as above. Assign each source generator a target expression of the corresponding type:
\[
    \lambda_t:S_{1,t}\longrightarrow(FA_{L_2})_{\alpha(t)},
    \qquad t\in T_1.
\]
Send each source basepoint generator to the corresponding target basepoint. Require every relation in $R_1$ to hold in $FA_{L_2}$ after translating the operation symbols by $j$ and the generators by $\lambda$. The universal property of the free construction, followed by the quotient by $R_1$, then gives a unique basepoint-preserving morphism
\[
    (j,\Phi):(\mathcal O_{L_1},FA_{L_1})
    \longrightarrow(\mathcal O_{L_2},FA_{L_2}).
\]

For $i=1,2$, let $\iota_i$ and $\operatorname{res}_i$ be the canonical embedding and restriction from that proposition:
\[
    \iota_i:A_{L_i}\hookrightarrow FA_{L_i},
    \qquad
    \operatorname{res}_i:FA_{L_i}\twoheadrightarrow A_{L_i},
    \qquad
    \operatorname{res}_i\circ\iota_i=\mathrm{id}_{A_{L_i}}.
\]
First restrict $\Phi$ to the source semantic subalgebra via $\iota_1$, then apply the target semantic restriction. This defines
\[
    e_t:=\operatorname{res}_{2,\alpha(t)}
        \circ\Phi_t\circ\iota_{1,t}:
    (A_{L_1})_t\longrightarrow(A_{L_2})_{\alpha(t)}.
\]

To obtain injective value correspondences, choose injections between the ordinary value domains
\[
    v_t:(A_{L_1})_t\setminus\{\ast\}
    \hookrightarrow
    (A_{L_2})_{\alpha(t)}\setminus\{\ast\},
\]
and prescribe the generator assignment on each ordinary source value symbol $\widehat a$ by
\[
    \lambda_t(\widehat a)
    =\iota_{2,\alpha(t)}(v_t(a)),
    \qquad a\in(A_{L_1})_t\setminus\{\ast\}.
\]
The remaining generating variables may be assigned well-typed target expressions, subject to preservation of $R_1$.

\begin{lemma}
\label{lem:semantic-interpretation-from-free-algebras}
Let $(j,\Phi):(\mathcal O_{L_1},FA_{L_1})\to(\mathcal O_{L_2},FA_{L_2})$ be a basepoint-preserving morphism with type map $\alpha:T_1\hookrightarrow T_2$. Suppose that, for each $t\in T_1$, there is an injection
\[
    v_t:(A_{L_1})_t\setminus\{\ast\}
    \hookrightarrow
    (A_{L_2})_{\alpha(t)}\setminus\{\ast\}
\]
satisfying
\[
    \Phi_t(\iota_{1,t}(a))
    =\iota_{2,\alpha(t)}(v_t(a)),
    \qquad a\in(A_{L_1})_t\setminus\{\ast\}.
\]
Then the maps
\[
    e_t:=\operatorname{res}_{2,\alpha(t)}
        \circ\Phi_t\circ\iota_{1,t}
\]
define a semantically faithful morphism $(j,e)$ of semantic algebras.
\end{lemma}

\begin{proof}
Each factor in the definition of $e$ is a basepoint-preserving operadic algebra morphism. Their composite is therefore a morphism over $j$, and $e_t(\ast)=\ast$.

For an ordinary source value $a$, the assumed identity gives
\[
\begin{aligned}
    e_t(a)
    &=\operatorname{res}_{2,\alpha(t)}
      \bigl(\Phi_t(\iota_{1,t}(a))\bigr)\\
    &=\operatorname{res}_{2,\alpha(t)}
      \bigl(\iota_{2,\alpha(t)}(v_t(a))\bigr)\\
    &=v_t(a)\neq\ast.
\end{aligned}
\]
If $e_t(a)=e_t(b)$, then either both $a,b$ are the basepoint, or both are ordinary values. In the latter case, $v_t(a)=v_t(b)$ implies $a=b$ by injectivity of $v_t$. Hence every $e_t$ is injective.
\end{proof}

\begin{theorem}[Faithful Compilation]
\label{thm:faithful-compiler-construction}
For $i=1,2$, let $FA_{L_i}=F_{L_i}/{\sim_{R_i}}$ be the value-generator presentations of Proposition~\ref{prop:semantic-restriction-and-embedding}, with generator families $S_i$. Let $j:\mathcal O_{L_1}\hookrightarrow\mathcal O_{L_2}$ be an operad embedding with injective type map $\alpha:T_1\hookrightarrow T_2$, and choose maps
\[
    \lambda_t:S_{1,t}\longrightarrow(FA_{L_2})_{\alpha(t)}.
\]
Write $\lambda^\#$ for the extension to the raw free algebra $F_{L_1}$ obtained by translating basic operation symbols through $j$. Assume the following conditions:
\begin{enumerate}
    \item For every $t\in T_1$, the basepoint generator satisfies $\lambda_t(\widehat\ast_t)=\ast$.
    \item For every relation $u=v$ in $R_1$ of type $t$,
    \[
        \lambda^\#_t(u)=\lambda^\#_t(v)
        \quad\text{in }(FA_{L_2})_{\alpha(t)}.
    \]
    \item There is a family of injections
    \[
        v_t:(A_{L_1})_t\setminus\{\ast\}
        \hookrightarrow
        (A_{L_2})_{\alpha(t)}\setminus\{\ast\},
        \qquad t\in T_1,
    \]
    such that, for every $t\in T_1$ and $a\in(A_{L_1})_t\setminus\{\ast\}$,
    \[
        \lambda_t(\widehat a)
        =\iota_{2,\alpha(t)}(v_t(a)).
    \]
    \item The family with value set $v_t((A_{L_1})_t\setminus\{\ast\})\cup\{\ast\}$ at each type $\alpha(t)$ and $\{\ast\}$ at every other target type is closed under all basic target operations, including nullary operations.
\end{enumerate}
Then $\lambda$ induces a unique basepoint-preserving morphism $(j,\Phi)$ between the free semantic algebras. Set
\[
    e_t:=\operatorname{res}_{2,\alpha(t)}
        \circ\Phi_t\circ\iota_{1,t}.
\]
Then the extension by zero $j_!A_{L_1}$ is defined and the associated compiler
\[
    (j,c):(\mathcal O_{L_1},A_{L_1})
    \longrightarrow
    (\mathcal O_{L_2},j_!A_{L_1})
\]
is semantically faithful, and each component $c_t:(A_{L_1})_t\to(j_!A_{L_1})_{\alpha(t)}$ is a bijection. After restriction along $j$, these maps give an isomorphism $c:A_{L_1}\xrightarrow{\sim}j^*j_!A_{L_1}$ of $\mathcal O_{L_1}$-algebras.
\end{theorem}

\begin{proof}
Let $z\in(FA_{L_1})_t$ and choose a representative $u\in(F_{L_1})_t$. Set $\Phi_t(z)=\lambda^\#_t(u)$. By condition~(2), the two sides of each relation in $R_1$ have the same image under $\lambda^\#$. Since $\lambda^\#$ preserves the operations over $j$, replacing one side of a relation by the other in any expression preserves its image. Thus any two representatives of $z$ have the same image, so $\Phi$ is well-defined on the quotient algebra. It preserves operations because $\lambda^\#$ does, and it preserves basepoints by condition~(1). Its values on the generators are prescribed by $\lambda$, which determines $\Phi$ uniquely.

For an ordinary source value $a$, the generator $\widehat a$ represents $\iota_{1,t}(a)$. The definition of $\Phi$ and condition~(3) therefore give
\[
    \Phi_t(\iota_{1,t}(a))
    =\lambda_t(\widehat a)
    =\iota_{2,\alpha(t)}(v_t(a)).
\]
By Lemma~\ref{lem:semantic-interpretation-from-free-algebras}, $(j,e)$ is semantically faithful, with $e_t(a)=v_t(a)$ on ordinary values and $e_t(\ast)=\ast$. The family in condition~(4) thus consists of the images $e_t((A_{L_1})_t)$ at the translated types and $\{\ast\}$ elsewhere. It is closed under the basic target operations by condition~(4), and hence under their operadic composites. The structure maps of $A_{L_2}$ therefore restrict to this family and define $j_!A_{L_1}$.

At each source type $t$, pulling back along $j$ gives
\[
    (j^*j_!A_{L_1})_t
    =(j_!A_{L_1})_{\alpha(t)}
    =e_t((A_{L_1})_t)
    \cong (A_{L_1})_t ,
\]
where the last bijection sends $e_t(a)$ to $a$. For $o\in\mathcal O_{L_1}(t_1,\ldots,t_n;t)$ and $a_i\in(A_{L_1})_{t_i}$, the definitions of pullback and extension by zero, followed by compatibility of $e$ with $o$, give
\[
    \theta^{j^*j_!A_{L_1}}(o;e\mathbf a)
    =\theta^{A_{L_2}}(j(o);e\mathbf a)
    =e_t\bigl(\theta^{A_{L_1}}(o;\mathbf a)\bigr),
\]
where $e\mathbf a=(e_{t_1}(a_1),\ldots,e_{t_n}(a_n))$. Thus the bijections $e_t(a)\mapsto a$ give an isomorphism of $\mathcal O_{L_1}$-algebras from $j^*j_!A_{L_1}$ to $A_{L_1}$. Its inverse has components $c_t(a)=e_t(a)$, uniquely determined by $\iota_{\alpha(t)}\circ c_t=e_t$ since $\iota$ is an inclusion. These are the components of the compiler of Definition~\ref{def:compiler}, so the compiler is semantically faithful.
\end{proof}

The compiler indeed enables us to translate a program from $L_1$ to $L_2$.

Let $j_1:\mathcal O_X\hookrightarrow\mathcal O_{L_1}$ embed the types and operations of a program $X$ in $L_1$. We first extend $A_X$ by zero along $j_1$ to obtain the $\mathcal O_{L_1}$-algebra $j_{1!}A_X\subseteq A_{L_1}$.

The compiler $(j,c)$ of Theorem~\ref{thm:faithful-compiler-construction} induces the embedding $j_2=j\circ j_1:\mathcal O_X\hookrightarrow\mathcal O_{L_2}$. To give the compiled image $c(j_{1!}A_X)$ an $\mathcal O_{L_2}$-algebra structure, we require the compatibility conditions for extension by zero given before Definition~\ref{def:compiler}. When they hold, the resulting algebra is precisely $j_!(j_{1!}A_X)$.

\begin{proposition}
\label{prop:compiled-program-semantic-equivalence}
Suppose $X$ is faithfully represented in $L_1$, so that the source interpretation gives $A_X\cong j_1^*j_{1!}A_X$. If $j_{1!}A_X$ admits extension by zero along the compiler, restricting the resulting algebra in $L_2$ to $\mathcal O_X$ recovers the original semantic algebra through the canonical isomorphisms
\[
    A_X \cong j_2^*\bigl(j_!(j_{1!}A_X)\bigr) \cong j_2^*j_{2!}A_X .
\]
\end{proposition}

\begin{proof}
Since $j_2=j\circ j_1$, functoriality of extension by zero gives
\[
    j_!(j_{1!}A_X)\cong j_{2!}A_X .
\]
Pulling back along $j_2$ and using $j_2^*=j_1^*j^*$ gives
\[
    j_2^*j_{2!}A_X
    \cong j_1^*j^*j_!(j_{1!}A_X)
    \cong j_1^*j_{1!}A_X
    \cong A_X .
\]
The second isomorphism uses the unit isomorphism of Remark~\ref{rem:extension-by-zero} at $j_{1!}A_X$, and the last uses the faithful source representation.
\end{proof}

\section{Extending a Compiler by Finite Products and Recursion}
\label{subsec:compiler-extension}

The variety of types in \Cxx{} and \Lean{} makes it impractical to verify the conditions of Theorem~\ref{thm:faithful-compiler-construction} separately for each type. Once these conditions have been verified for the fundamental types and their basic operations, how can we extend the compiler to new types and operations while preserving semantic faithfulness?

We focus on finite Cartesian products and recursion, which together suffice for most of the applications considered here. For example, a fixed-length vector is a product of copies of its element type, and the sum of its entries is defined recursively from zero and addition. For a product, the translations of the components determine the translation of a tuple. For recursion, it suffices to check that translation preserves the initial operation and the step.

Let $T$ be the set of colors corresponding to the fundamental types of a programming language, let $\mathcal O$ be the $T$-colored operad generated by the basic operations on these types, and let $A=(A_t)_{t\in T}$ be the family of their ordinary value sets. For any finite list of colors $t_1,\ldots,t_n$, allowing repetitions and $n=0$, take the product color $p=\prod_{i=1}^n t_i$ with ordinary value set $A_p=\prod_{i=1}^n A_{t_i}$. For $n=1$ we use $t_1$ itself, with the identity as both constructor and projection. Extend $\mathcal O$ by a tuple constructor $\mathsf{tuple}:(t_1,\ldots,t_n)\to p$ and projections $\pi_i:p\to t_i$. On ordinary values they satisfy the equations
\[
    \pi_i\bigl(\mathsf{tuple}(a_1,\ldots,a_n)\bigr)=a_i,
    \qquad
    \mathsf{tuple}\bigl(\pi_1(z),\ldots,\pi_n(z)\bigr)=z.
\]
For any set $S$, maps $h_i:S\to A_{t_i}$ therefore determine a unique map $h:S\to A_p$ with $\pi_i\circ h=h_i$, given by $h(s)=\mathsf{tuple}(h_1(s),\ldots,h_n(s))$. In particular, the uniqueness part of this universal property gives the equality criterion for product values:
\[
    z=z' \quad\Longleftrightarrow\quad
    \pi_i(z)=\pi_i(z')\ \text{for every }i.
\]
For $n=0$ the indexing diagram is empty. Every set $S$ admits exactly one cone over it, so its limit is terminal in $\mathbf{Set}$. Writing $\mathbf{1}$ for the empty product color, its ordinary carrier is $A_{\mathbf{1}}=\{()\}$, and the unique map $S\to A_{\mathbf{1}}$ sends each element to the empty tuple $()$. For example, fixed-length vectors of color $t$ can be formed recursively by $V_0=\mathbf{1}$ and $V_{n+1}=t\times V_n$.

Choosing ordinary values $d_j\in A_{t_j}$ for $j\ne i$ also lets us adjoin a coordinate embedding $\operatorname{ins}_i:t_i\to p$. It places its argument in coordinate $i$ and $d_j$ in each other coordinate, so $\pi_i(\operatorname{ins}_i(a))=a$ for every $a\in A_{t_i}$. As in Definition~\ref{def:semantic-algebra}, we adjoin the undefined outcome $\ast$ to $A_p$ and take the strict totalizations of the new operations.

\begin{proposition}
\label{prop:product-type-algebra}
Let $A^\ast$ be the given $\mathcal O$-algebra on the fundamental types. For a finite product $p=\prod_{i=1}^n t_i$, let $\mathcal O'$ be the resulting operad on $T'=T\cup\{p\}$. For $n=0$ adjoin just the nullary tuple constructor, for $n=1$ use $\mathcal O$ itself, and for $n\geq 2$ freely adjoin the tuple constructor, projections and chosen coordinate embeddings subject to $\pi_i\circ\operatorname{ins}_i=\mathbf 1_{t_i}$. With
\[
    A_p^\ast=A_p\cup\{\ast\}=\Bigl(\prod_{i=1}^n A_{t_i}\Bigr)\cup\{\ast\},
\]
the extended family $(A_t^\ast)_{t\in T'}$ is an $\mathcal O'$-algebra whose restriction along $\mathcal O\to\mathcal O'$ is $A^\ast$.
\end{proposition}
\begin{proof}
Interpret the old operations by the structure maps of $A^\ast$ and the new basic operations by their strict totalizations. For the empty product the nullary constructor takes the value $()$. By Proposition~\ref{prop:semantic-algebra}, these interpretations respect substitution, units, and input permutations on the freely generated operad.

The relations among the old operations hold in $A^\ast$. For each chosen coordinate embedding, the new relation holds on values since
\[
    \pi_i(\operatorname{ins}_i(z))=z \qquad (z\in A_{t_i}^\ast).
\]
For an ordinary $z$ the projection recovers the inserted coordinate, and for $z=\ast$ both maps preserve the basepoint. The structure maps therefore descend to the quotient operad $\mathcal O'$ and define the asserted algebra. Their restriction to $\mathcal O$ is the given algebra structure on $A^\ast$. The same construction can now be applied to products of the extended types.
\end{proof}

For the extensions below, fix a compiler $(j,c)$ satisfying the hypotheses of Theorem~\ref{thm:faithful-compiler-construction}, with type map $\alpha$. We write $e_t=\iota_{\alpha(t)}\circ c_t$ for the interpretation of source values in $A_{L_2}$. Its operad embedding $j$ is induced by an embedding of the basic signatures. We use primes for target data and match each newly adjoined color and symbol with a fresh target counterpart. On each new color we extend the generator assignment by $\widehat a\mapsto\widehat{e_t(a)}$ for values, including the basepoint, and choose corresponding target variables for the new variables. Any further basic operation freely adjoined on the target side must preserve the extended image family, including in arity zero.

\begin{proposition}
\label{prop:compiler-product-extension}
For a finite list of types $t_1,\ldots,t_n$ of $L_1$, extend the source semantic algebra $A_{L_1}$ by the product type $p=\prod_i t_i$ and the target semantic algebra $A_{L_2}$ by $p'=\prod_i\alpha(t_i)$ as in Proposition~\ref{prop:product-type-algebra}. For each chosen coordinate embedding, take the target filler corresponding to $d_j$ to be $e_{t_j}(d_j)$. Then the compiler $(j,c)$ from $L_1$ to $L_2$ extends to a semantically faithful compiler between the extended semantic algebras. It sends $p$ to $p'$, and its value map $e_p$ is uniquely determined by the component maps $e_{t_i}$.
\end{proposition}
\begin{proof}
On ordinary values, preservation of the tuple constructor forces
\[
    e_p(a_1,\ldots,a_n)=(e_{t_1}(a_1),\ldots,e_{t_n}(a_n)).
\]
Preservation of the basepoint gives $e_p(\ast)=\ast$. For the empty product these two requirements give $e_p(())=()$ and $e_p(\ast)=\ast$, respectively. For a unary product they recover the given value map. Since each $e_{t_i}$ is injective and sends ordinary values to ordinary values, $e_p$ is injective.

Taking the $i$th coordinate gives $\pi'_i\circ e_p=e_{t_i}\circ\pi_i$. The choice of the other coordinates gives $e_p\circ\operatorname{ins}_i=\operatorname{ins}'_i\circ e_{t_i}$ for each chosen embedding. If an input is $\ast$, the source operation and its target counterpart both return $\ast$ by strictness. Thus the new concrete computation relations and basepoint propagation relations are preserved.

Since the correspondence of product generators preserves and reflects the defining cancellations, it extends $j$ to an embedding of the extended operads. The image on the product color is
\[
    \operatorname{im}(e_p)=\{\ast\}\,\cup\,\prod_{i=1}^n\bigl(\operatorname{im}(e_{t_i})\setminus\{\ast\}\bigr).
\]
The target tuple constructor takes inputs from the component images into $\operatorname{im}(e_p)$, and each projection takes $\operatorname{im}(e_p)$ into the corresponding component image. The chosen fillers give closure under the coordinate embeddings as well. Thus Theorem~\ref{thm:faithful-compiler-construction} gives a semantically faithful compiler on the extended algebras. Its restriction to the original types and operations is $(j,c)$.
\end{proof}

For the recursion construction, let $I$, $P$, and $S$ be colors for an index, a parameter, and a state, respectively. A step $s:(I,P,S)\to S$ takes the current index $k$, the parameter $x$, and the current state $a$ to the next state. By Proposition~\ref{prop:product-type-algebra}, we can pack its inputs into the product color $I\times P\times S$. The induced function on product values sends an ordinary tuple $(k,x,a)$ to $s(k,x,a)$ and the basepoint to the basepoint.

The colors $P$ and $S$ may coincide. For example, take $I=\texttt{Nat}$ and $P=S=\texttt{Int}$ in \Lean{}. Starting at an integer $x$ and incrementing the state at each step uses the initial operation $b(x)=x$ and the step $s(k,x,a)=a+1$. Here $x$ supplies the initial state and $a$ is the value updated at each step.

\begin{proposition}
\label{prop:parameterized-recursion}
Suppose the index color $I$ has ordinary value domain $\mathbb N$. Extend the operad on its existing colors by an initial operation $b:P\to S$ and a step operation $s:(I,P,S)\to S$. For every choice of their strict interpretations, there is a unique map $r:A_I^\ast\times A_P^\ast\to A_S^\ast$ that is strict in each argument and satisfies, for every $k\in\mathbb N$ and $x\in A_P^\ast$,
\[
    r(0,x)=b(x),\qquad r(k+1,x)=s\bigl(k,x,r(k,x)\bigr).
\]
Here $b$ and $s$ also denote their interpretations. This universal property characterizes $r$ from the two specified operations $b$ and $s$.
\end{proposition}
\begin{proof}
Fix $x$ and construct a sequence of states by taking $a_0=b(x)$ and $a_{k+1}=s(k,x,a_k)$. Set $r(k,x)=a_k$ and $r(\ast,x)=\ast$. When $x=\ast$, strictness of $b$ gives $a_0=\ast$ and strictness of $s$ gives $a_{k+1}=\ast$ at every step. Thus $r$ is strict and satisfies both equations. For any other strict solution, the initial equation fixes the same state at zero. If the states agree at $k$, applying $s$ to the same index, parameter, and state gives agreement at $k+1$. Induction gives equality at every ordinary index, and strictness gives equality at $\ast$.
\end{proof}

We can also use a universal property to specify the index type itself. If the natural numbers are to be introduced as a new color $N$, adjoin constructors $0:()\to N$ and $\mathsf{succ}:N\to N$. Require the ordinary carrier $A_N$ to have the following property. For every set $Y$, every $y_0\in Y$, and every map $u:Y\to Y$, there is a unique $h:A_N\to Y$ such that
\[
    h(0)=y_0,\qquad h(\mathsf{succ}(n))=u(h(n)).
\]
For the ordinary natural numbers, applying $u$ to $y_0$ exactly $n$ times defines $h(n)$, and induction proves uniqueness. The property determines the type with its constructors up to unique isomorphism. Indeed, applying it in both directions between two such types gives constructor-preserving maps whose composites are identities by uniqueness. In our pointed semantics we then adjoin $\ast$ to the ordinary carrier and extend $\mathsf{succ}$ strictly.

To recover Proposition~\ref{prop:parameterized-recursion} from this property, fix $x$ and take the set $Y=\mathbb N\times A_S^\ast$, with $y_0=(0,b(x))$ and $u(k,a)=(k+1,s(k,x,a))$. The unique map $h$ then satisfies $h(k)=(k,r(k,x))$ by induction. Its second projection gives the required state after $k$ steps.

This is the constructor pattern of \Lean{}'s \texttt{Nat.zero} and \texttt{Nat.succ}. Its recursion principle defines operations on that type, including addition by recursion on the second argument. \Lean{} also introduces \texttt{Int} inductively, using two constructors $\texttt{Int.ofNat},\texttt{Int.negSucc}:\texttt{Nat}\to\texttt{Int}$ for the integers $n$ and $-(n+1)$, respectively~\cite[chapters on natural numbers and integers]{lean-reference}.

\begin{proposition}
\label{prop:compiler-recursion-extension}
Extend the source semantic algebra $A_{L_1}$ by the initial and step operations $b,s$ and the target semantic algebra $A_{L_2}$ by $b',s'$ as in Proposition~\ref{prop:parameterized-recursion}. If the natural-number index colors $I,I'$ are already present, assume that $e_I$ preserves zero and successor. If they are newly adjoined, take $e_I$ to be the unique map preserving zero, successor and the basepoint. Suppose that
\[
\begin{aligned}
    e_S(b(x))     &= b'(e_P(x)),\\
    e_S(s(k,x,a)) &= s'(e_I(k),e_P(x),e_S(a)).
\end{aligned}
\]
Then the compiler $(j,c)$ from $L_1$ to $L_2$ extends to a semantically faithful compiler between the extended semantic algebras. For the recursive operations $r$ and $r'$ determined by $(b,s)$ and $(b',s')$, respectively, this extension satisfies
\[
    e_S(r(k,x))=r'(e_I(k),e_P(x)).
\]
\end{proposition}
\begin{proof}
For newly adjoined index colors, the universal property gives a constructor-preserving map on ordinary values in each direction. Uniqueness makes their composites the identity maps. Extending this bijection by $e_I(\ast)=\ast$ gives the asserted value map, whose image contains zero and is closed under successor.

Fix an ordinary parameter $x$. The recursion equations for $r$ and compatibility with $b$ and $s$ give
\[
\begin{aligned}
    e_S(r(0,x))
        &= e_S(b(x)) = b'(e_P(x)),\\
    e_S(r(k+1,x))
        &= e_S\bigl(s(k,x,r(k,x))\bigr)
         = s'\bigl(e_I(k),e_P(x),e_S(r(k,x))\bigr).
\end{aligned}
\]
Since $e_I$ preserves zero and successor, the recursion equations for $r'$ give the same initial value and step for $k\mapsto r'(e_I(k),e_P(x))$. Uniqueness in Proposition~\ref{prop:parameterized-recursion} therefore gives $e_S(r(k,x))=r'(e_I(k),e_P(x))$ for every $k\in\mathbb N$. If either input is $\ast$, both sides equal $\ast$ by strictness and preservation of basepoints.

When the initial and step operations are corresponding composites of previously verified operations, their compatibility follows from preservation of substitution. Matching the freely adjoined operation symbols extends the operad embedding. If the signatures also include symbols for $r$ and $r'$, we interpret them by the unique maps of Proposition~\ref{prop:parameterized-recursion}. The compatibility identities for the new operations preserve their concrete computation relations and place each output on translated inputs in the corresponding image. Strictness gives the basepoint propagation relations. Theorem~\ref{thm:faithful-compiler-construction} now gives the semantically faithful extension of $(j,c)$.
\end{proof}

\begin{theorem}
\label{thm:compiler-extension-products-recursion}
A compiler satisfying the hypotheses of Theorem~\ref{thm:faithful-compiler-construction} on the fundamental types extends to a semantically faithful compiler after any finite sequence of the product and recursion extensions of Propositions~\ref{prop:compiler-product-extension} and~\ref{prop:compiler-recursion-extension}.
\end{theorem}
\begin{proof}
Induct on the number of extensions. At each step the corresponding proposition carries the hypotheses of Theorem~\ref{thm:faithful-compiler-construction} to the extended algebras, so the resulting compiler is semantically faithful.
\end{proof}

\section{Local-to-Global Reconstruction of Program Semantics}
\label{sec:local-to-global}

How does the complier translate a program and which data structure is needed to represent the program? If we use the Theorem~\ref{thm:faithful-compiler-construction} directly, we need to expand all the informtaion that the program contains, for example, if the program calls a function twice, we need to expand the function twice. This is not efficient, all we need is to record the function once and record how the program calls the function. We will discuss how to use the idea in sheaf theory to reduce the work of verifying the conditions of Theorem~\ref{thm:faithful-compiler-construction} and how to represent a program in a data structure that is suitable for compilation.

In mathematics and physics, we have defined such an object, called a factorization algebra, which is a cosheaf that allows us to assemble global data from local data~\cite[Definitions~1--3]{costello-gwilliam-factorization}.
In algebraic topology, Ayala and Francis use factorization homology to prove a common generalization of Poincar\'e duality for manifolds and Koszul duality for \(E_n\)-algebras~\cite{ayala-francis-poincare-koszul}.
In four-dimensional Chern--Simons theory with simple gauge Lie algebra \(\mathfrak g\), Costello uses Koszul duality to recover Drinfeld's Yangian \(Y(\mathfrak g)\) from the algebra of local observables~\cite[Thm.~9.5.1]{costello-yangian}.
In algebraic quantum field theory on globally hyperbolic spacetimes, Benini, Perin, and Schenkel recover the local associative products from time-orderable prefactorization algebras satisfying additivity and the time-slice axiom~\cite[Thm.~5.1]{benini-perin-schenkel-aqft}.

In the context of programming languages, we can use a similar idea to represent a program as a factorization algebra, where the local data corresponds to the pragment of the program and the global data corresponds to the entire program. This allows us to efficiently represent and compile programs by only recording the necessary information of the local data and how they are assembled to form the global data. In this section, we will define a factorization algebra for programming languages and show how it can be used to represent and compile programs.

\begin{definition}[Prefactorization algebra]
\label{def:prefactorization-algebra}
Let $M$ be a topological space and $(\mathcal C,\otimes,\mathbf 1)$ a symmetric monoidal category. A \emph{prefactorization algebra} $\mathcal F$ assigns an object $\mathcal F(U)$ to each open $U\subseteq M$. For every finite family of pairwise disjoint opens $U_1,\ldots,U_n\subseteq V$, it gives a morphism
\[
    \mu^V_{U_1,\ldots,U_n}:
    \mathcal F(U_1)\otimes\cdots\otimes\mathcal F(U_n)
    \longrightarrow \mathcal F(V).
\]
Reordering the opens permutes the tensor factors by the corresponding symmetry. The unary map $\mu^U_U$ is the identity. If $V_1,\ldots,V_k$ are disjoint in $W$ and the $U_{ij}$ are disjoint in each $V_i$, then, up to the canonical tensor identifications,
\[
    \mu^W_{V_1,\ldots,V_k}\circ
    \bigotimes_{i=1}^k \mu^{V_i}_{U_{i1},\ldots,U_{in_i}}
    =\mu^W_{(U_{ij})_{i,j}} .
\]
In particular the unary maps form a covariant functor on opens. We allow the empty family, whose source is $\mathbf 1$. Costello and Gwilliam give this definition for vector spaces with their tensor product~\cite[Definition~1]{costello-gwilliam-factorization}.
\end{definition}

The original definition of factorization algebra is based on the notion of topological spaces, that is we need to give a topology on the space of the program. In our case, we can define the open set of the program as the sub-program of the program, and the open cover of the program as the collection of sub-programs that cover the entire program. The prefactorization algebra is then defined as a functor from the category of open sets of the program to the category of semantic algebras.

Use the program $X$ itself as the underlying space, with its syntax locations as points. An open $U\subseteq X$ contains the dependencies of every location it contains. Unions and intersections preserve this condition. To model separate function bodies, distinguish a function's declared interface from its body. A call depends on that interface and on its argument expressions. Its own fragment can therefore refer to the function without containing the body. The fragment for the definition also contains the interface.

Fix the semantic presentation $(\mathcal O_X,A_X)$ of Remark~\ref{rem:program-semantic-algebras}. At each location $x\in X$, let $S_{x,t}$ consist of the values of type $t$ supplied at $x$, such as a literal or the admissible inputs at a parameter. For each open $U\subseteq X$, put
\[
    S_{U,t}=\bigcup_{x\in U}S_{x,t}.
\]
Write $\langle S_U\rangle_X$ for the pointed $\mathcal O_X$-subalgebra of $A_X$ generated by $S_U$. We assume that the values supplied throughout the program generate $A_X$, so $\langle S_X\rangle_X=A_X$.

Let $T_U\subseteq T_X$ contain the types occurring in $U$ and every type at which $\langle S_U\rangle_X$ has an ordinary value. The basic operation symbols occurring in $U$, with their original input and output types, determine an inclusion $j_U:\mathcal O_U\hookrightarrow\mathcal O_X$. Restricting along $j_U$ and extending by zero inside $A_X$ gives
\begin{equation}
    \label{eq:fragment-generated-carriers}
    A_U=j_U^*\langle S_U\rangle_X,
    \qquad
    \mathcal F_X(U):=j_{U!}A_U=\langle S_U\rangle_X.
\end{equation}
By the unit isomorphism $j_U^*j_{U!}\cong\mathrm{id}$ of Remark~\ref{rem:extension-by-zero}, we have $j_U^*j_{U!}A_U\cong A_U=j_U^*\langle S_U\rangle_X$ via the given inclusions. Outside $T_U$, both have value set $\{\ast\}$, so they agree as subalgebras of $A_X$.

Extension by zero uses the inclusions $(A_U)_t\hookrightarrow(A_X)_t$ for $t\in T_U$. Thus the local data on $U$ record how the fragment sits in the program through the embedding $\mathcal F_X(U)=j_{U!}A_U\hookrightarrow A_X$. We assemble these algebras in the category $\operatorname{Sub}_{\mathcal O_X}(A_X)$ of pointed $\mathcal O_X$-subalgebras and inclusions. Its coproduct $B\vee C=\langle B\cup C\rangle_X$ is the least subalgebra containing both, with the union taken type by type. We use this coproduct as the tensor product. Its unit is $\mathcal F_X(\varnothing)=\langle\varnothing\rangle_X$, the pointed subalgebra generated by the interpretations of ground terms.

\begin{proposition}[]
\label{prop:program-prefactorization}
The extensions by zero $\mathcal F_X(U)$ of~\eqref{eq:fragment-generated-carriers} define a prefactorization algebra on $X$ valued in $\operatorname{Sub}_{\mathcal O_X}(A_X)$.
\end{proposition}
\begin{proof}
If $U\subseteq V$, then $S_U\subseteq S_V$, so $\mathcal F_X(U)\subseteq\mathcal F_X(V)$. This inclusion preserves every $\mathcal O_X$-operation because the two interpretations restrict the same operation of $A_X$.

For pairwise disjoint opens $U_1,\ldots,U_n\subseteq V$, the subalgebra $\mathcal F_X(V)$ thus contains each $\mathcal F_X(U_j)$. Taking the subalgebra generated by their union gives
\[
    \bigvee_j\mathcal F_X(U_j)
    =\Big\langle\bigcup_j\mathcal F_X(U_j)\Big\rangle_X
    \subseteq\mathcal F_X(V).
\]
Since the tensor product is this coproduct, the displayed inclusion is the unique morphism whose restriction to each $\mathcal F_X(U_j)$ is its inclusion into $\mathcal F_X(V)$. We take it as the structure map $\mu^V_{U_1,\ldots,U_n}$. Reordering the opens leaves their join unchanged. Successive inclusions compose to the inclusion into the final fragment, giving the substitution and unit identities.
\end{proof}

\begin{definition}[Factorization algebra]
\label{def:factorization-algebra}
A family $\{V_i\}_{i\in I}$ of opens in $U$ is a \emph{Weiss cover} if every finite subset of $U$ lies in some $V_i$. Following the underived definition of Costello and Gwilliam~\cite[Definitions~2 and~3]{costello-gwilliam-factorization}, a prefactorization algebra valued in a category with small coproducts is a \emph{factorization algebra} if it satisfies two conditions.
\begin{enumerate}
    \item For disjoint opens $U_1,U_2$, the structure map
    $\mathcal F(U_1)\otimes\mathcal F(U_2)\to\mathcal F(U_1\sqcup U_2)$
    is an isomorphism.
    \item For every Weiss cover of $U$, the two maps induced by the overlaps have $\mathcal F(U)$ as their coequalizer in $\mathcal C$.
    \[
        \coprod_{i,j\in I}\mathcal F(V_i\cap V_j)
        \;\rightrightarrows\;
        \coprod_{i\in I}\mathcal F(V_i)
        \longrightarrow \mathcal F(U).
    \]
    The parallel arrows use the inclusions into $V_i$ and $V_j$. Thus maps from the $\mathcal F(V_i)$ to a common object that agree on all overlaps induce a unique map from $\mathcal F(U)$. In vector spaces this says that the direct sum of the local spaces, modulo the differences of the two images from each overlap, is $\mathcal F(U)$.
\end{enumerate}
\end{definition}

The gluing data is then given by the morphisms of semantic algebras that relate the local data to the global data. The factorization algebra is then defined as a prefactorization algebra that satisfies the gluing condition, which ensures that the local data can be assembled to form the global data. 

\begin{proposition}[]
\label{prop:program-factorization}
With the inclusion maps and tensor product above, $\mathcal F_X$ is a factorization algebra with global value $\mathcal F_X(X)=A_X$. For every ordinary open cover $U=\bigcup_i V_i$,
\[
    \mathcal F_X(U)=\bigvee_i\mathcal F_X(V_i)
    \quad\text{in }\operatorname{Sub}_{\mathcal O_X}(A_X).
\]
In particular, taking the functions of $X$ as the open sets $U_f$ and the call relations between them as the transition maps, we obtain
\[
    A_X=\mathcal F_X(X)=\bigvee_f\mathcal F_X(U_f).
\]
\end{proposition}
\begin{proof}
Let $U_1$ and $U_2$ be disjoint opens and put $U=U_1\sqcup U_2$. Since the tensor product is the join, the structure map of Proposition~\ref{prop:program-prefactorization} is the inclusion
\[
    \mathcal F_X(U_1)\vee\mathcal F_X(U_2)\hookrightarrow\mathcal F_X(U).
\]
Conversely, the join contains both $S_{U_1}$ and $S_{U_2}$, hence their union $S_U$. As a subalgebra of $A_X$, it therefore contains $\langle S_U\rangle_X=\mathcal F_X(U)$. Thus the structure map is an isomorphism, as required by condition~(1) of Definition~\ref{def:factorization-algebra}.

For an arbitrary open cover $U=\bigcup_i V_i$, we have $S_U=\bigcup_i S_{V_i}$. Hence
\[
    \bigvee_i\mathcal F_X(V_i)
    =\Big\langle\bigcup_i S_{V_i}\Big\rangle_X
    =\langle S_U\rangle_X
    =\mathcal F_X(U).
\]
Inclusions $\mathcal F_X(V_i)\hookrightarrow C$ agree on the overlaps, and the join identity gives the unique inclusion $\mathcal F_X(U)\hookrightarrow C$ required by condition~(2). Since $S_X$ generates $A_X$, we have $\mathcal F_X(X)=\langle S_X\rangle_X=A_X$.

In the function case, every location of $X$ belongs to the definition of some $f$ or to its dependencies. Hence $X=\bigcup_f U_f$, and taking $U=X$ in the join formula gives the last equality.
\end{proof}

\begin{example}
Suppose a function $f$ has a formal parameter named \texttt{x}, and a caller with its own variable \texttt{x} makes the call \texttt{f(x)}. Write $(U_1,\mathtt{x})$ for the parameter and $(U_2,\mathtt{x})$ for the caller's variable. Lexical scope resolves each name to its declaration in the source program. Under the embeddings $\mathcal F_X(U_i)=j_{U_i!}A_{U_i}\hookrightarrow A_X$, the two bindings retain their distinct declaration locations $d_1$ and $d_2$, even when they have the same value. If both fragments refer to a global variable declared at $d$, then $d\in U_1\cap U_2$. The two embeddings agree on this overlap, so both references denote the same global binding after assembly.
\end{example}

\section{Implementing a Local-to-Global Compiler}
\label{sec:implementation}

We assume that a compiler for the fundamental types and their basic operations has already been constructed and satisfies the hypotheses of Theorem~\ref{thm:faithful-compiler-construction}. By Theorem~\ref{thm:compiler-extension-products-recursion}, this compiler extends to the compound types and operations obtained by any finite sequence of the product and recursion constructions of Section~\ref{subsec:compiler-extension}. The extension preserves the semantics at each step. To compile programs that use these types and operations, it remains to translate the computations within each function body and the calls between functions.

We base our compiler architecture on the factorization algebra $\mathcal F_X$ of Section~\ref{sec:local-to-global}. Each open in the function cover gives a unit of translation. We first describe the data attached to a fragment, then follow a variable through its lifetime and through the branches and loops of the function. At a call, the callee's declaration lies in the overlap of the two fragments. It identifies the function to invoke, and its parameter positions determine how to bind the caller's arguments. Section~\ref{subsec:compiler-workflow} shows how these local translations assemble into a \Lean{} program.

\subsection{Function Fragments}
\label{subsec:function-fragments}

For a function $f$ of the program $X$, let $U_f$ be the open set consisting of its definition and the dependencies of that definition. If $f$ calls $g$, the declaration of $g$ belongs to both $U_f$ and $U_g$. The call in $U_f$ refers to that declaration, while the body of $g$ is stored once in its own fragment. By Proposition~\ref{prop:program-factorization}, these opens cover $X$ and their local algebras generate the program algebra
\[
    A_X=\bigvee_f\mathcal F_X(U_f).
\]

The extractor stores one SQLite file per source file. Within it, each function definition has its own logical database, identified by \texttt{function\_id}. Its rows give the formal parameters in declaration order and the statements in source order within each region. A control statement links to its child regions. A call records the callee's declaration when resolved and the source range of each actual argument, indexed by parameter position. Thus two calls to $g$ in $f$ have separate argument records and refer to the same declaration of $g$.

\subsection{Variable Lifetimes in a Function Fragment}
\label{subsec:function-body-translation}

While translating a function body, the compiler keeps its active lexical scopes on a stack. It pushes a scope on entering a block and pops it on leaving the block. To resolve a local variable, it searches from the top for the nearest visible declaration of that name. Once resolved, a read or write in the generated code refers directly to the corresponding binding.

At run time, each invocation has its own frame. A local mutable binding $x_U$ of type $t$ owned by that invocation occupies one slot. Write $L=\mathbb N$ for the live states, with $0$ as the initial state. A live slot holds a value $a\in A_t$ and a state $\ell\in L$. An empty slot is denoted by $\ast$, so the lifetime state takes values in $L\cup\{\ast\}$. The transitions of a single slot are
\[
    \ast \xrightarrow{\mathrm{init}(a)} (a,0),
    \qquad (a,\ell) \xrightarrow{\mathrm{write}(b)} (b,\ell+1),
    \qquad (a,\ell) \xrightarrow{\mathrm{clear}} \ast .
\]
A read of a live slot returns $a$ and leaves $\ell$ unchanged. Clearing the slot ends the current lifetime. A subsequent initialization stores the supplied value and starts a new lifetime at $0$. Leaving a scope clears its owned local bindings, and returning from a function releases its frame.

For example, if the slot for \texttt{x} initially holds $(2,0)$, the statement \texttt{x = x + 1;} writes $(3,1)$ into that slot. A subsequent \texttt{return x;} reads $3$. If the assignment occurs inside a loop and \texttt{x} was declared outside it, each iteration updates the same binding. After $k$ writes in that lifetime its state is $k$.

The database records the position of each instruction within its region. A lifetime-end record identifies the declaration being cleared and the exit where its lifetime ends. The translator reads the instructions in order, so each use sees the latest value of its binding on the path taken. Consider \texttt{x = x + 1; y = x; x = x * 2;}. Write $x_0$ for the incoming value of \texttt{x}. The first write produces $x_1$ from $x_0$. The intervening read supplies $x_1$ to \texttt{y}, and the second write also uses $x_1$ to produce $x_2$.

\subsection{Control Flow within a Function Fragment}
\label{subsec:control-flow-translation}

Within $U_f$, we treat each branch of a conditional and each loop as a local function. For a branch, take its statements and the test that selects it. For a loop, take its test and body together with any initialization and update. These source regions are linked to the control statement in the database. Adjoining the dependencies of their locations gives an open set $W$. Since $U_f$ contains those locations and is closed under dependencies, $W\subseteq U_f$, with the induced inclusion $\mathcal F_X(W)\hookrightarrow\mathcal F_X(U_f)$. For example, if both branches read a variable declared before the conditional, that declaration belongs to both opens.

For an \texttt{if} statement, the compiler translates the two branches as local functions $h_{\mathrm{then}}$ and $h_{\mathrm{else}}$. The generated code evaluates the condition once. Both functions capture the bindings visible after the test, including any updates made by the condition. For a test result $d$, the \Lean{} operation \texttt{Cpp.chooseBranch} satisfies
\[
    \mathsf{chooseBranch}(d,h_{\mathrm{then}},h_{\mathrm{else}})=
    \begin{cases}
        h_{\mathrm{then}}(), & d=\mathrm{true},\\
        h_{\mathrm{else}}(), & d=\mathrm{false},\\
        \ast, & d=\ast.
    \end{cases}
\]
On normal completion, the selected branch returns the current values of the outer bindings that either branch may update. The statements following the conditional use these values. If the source omits the \texttt{else} clause, the false branch passes through the values left by the test. A \texttt{return} in the selected branch returns from the enclosing function. Thus corresponding test results select corresponding branches, and preservation of the selected branch's outcome gives preservation of the conditional.

We assume that every \texttt{while} or \texttt{for} loop reached on an admitted source input terminates. For each loop the compiler extracts a function $L_q$ with a termination test $q$, fixed parameters $p$, and a state $z$ carried from one iteration to the next. The parameters supply the values used by the test or body that remain fixed during the loop. The state contains the current mutable bindings they use, including any shared with the enclosing function. The compiler emits the body once and can specialize $L_q$ to the test at that loop site.

For \texttt{while (condition) body}, $q$ evaluates the source condition once and signals exit when the condition is false. Write $(d,z_g)=q(p,z)$, where $d$ is the exit flag and $z_g$ is the state after the test, including any writes performed by the condition. Let $s$ execute one iteration. When control proceeds normally through the body, the extracted function satisfies
\[
    L_q(p,z)=
    \begin{cases}
        z_g,                 & d=\mathrm{true},\\
        L_q(p,s(p,z_g)),     & d=\mathrm{false}.
    \end{cases}
\]
Only the selected branch executes. Thus an initially false source condition returns the state after the first test, and a successful iteration passes its updated state to the next test. An undefined operation on the executed path yields $\ast$.

For \texttt{for (init; condition; update) body}, the enclosing function executes \texttt{init} once to obtain the initial state. One full step $s$ executes \texttt{body} and then \texttt{update}. A \texttt{continue} directed at this loop goes to \texttt{update} in a \texttt{for} loop and directly to the next test in a \texttt{while} loop. A \texttt{break} exits the loop with the current state. A \texttt{return} exits the enclosing function according to Section~\ref{subsec:parameter-binding}. Both skip the remaining work of the iteration. An omitted \texttt{for} condition always permits another iteration, and an omitted update leaves the state unchanged. Bindings owned by a scope are cleared when control leaves it, including the loop index declared in \texttt{init} on exit from the loop~\cite[\texttt{[stmt.iter]}, \texttt{[stmt.jump]}]{cpp23}.

For example, consider \texttt{for (int i = 0; i < n; ++i) sum += i;} with $n\geq 0$. The fixed parameter is $n$. Writing just the values of the mutable slots, the initial state is $(0,u_0)$, where $u_0$ is the incoming value of \texttt{sum}. The test exits when $i\geq n$, and a full step sends $(i,u)$ to $(i+1,u+i)$ on defined inputs. Each write increments the lifetime counter of its own slot. If $n=0$, the enclosing function retains $u_0$. For $n=3$ and $u_0=0$, it resumes with \texttt{sum} equal to $3$ and the local slot for \texttt{i} cleared.

For a loop with a prescribed iteration count, as in this example, the initialization and full step give the operations $b$ and $s$ of Proposition~\ref{prop:compiler-recursion-extension}. For a data-dependent count, consider a defined source execution and induct on the number of body executions remaining. Assume that the translations of initialization, test, body and update preserve their results, control outcomes and observable slot states. Corresponding states then give the same exit flag and corresponding states after the test. An exit returns those states immediately. If the body executes, its translation gives the same control outcome and corresponding updates. A \texttt{break} or \texttt{return} finishes the comparison directly. Whenever control reaches the next test, one fewer body execution remains. Thus the target loop terminates with the corresponding outcome and slot states. The iteration count used in this proof may depend on the input.

\subsection{Parameter Binding at Function Calls}
\label{subsec:parameter-binding}

Let $U$ and $V$ be open fragments of a program $X$. Suppose a call in $U$ passes a variable named $x$ to a function defined in $V$, whose corresponding formal parameter is named $y$. Write $x_U$ for the binding of $x$ at its declaration in the caller's current invocation and $y_V$ for the formal at its declaration. These bindings remain distinct when the two declarations use the same spelling. The call interface records which actual argument supplies each formal parameter.

For a type $t$, let $N_t$ be the set of typed binding names and take the live states $L=\mathbb N$ from Section~\ref{subsec:function-body-translation}. Using the finite products of Proposition~\ref{prop:product-type-algebra}, represent a mutable variable and a constant value binding by the respective ordinary carriers
\[
    M_t=N_t\times A_t\times L,\qquad C_t=N_t\times A_t .
\]
Thus $(x_U,a,\ell)$ records the caller's binding, its current value $a$, and its lifetime state $\ell$. A read or write through a mutable binding requires a live state. Writing a value $b\in A_t$ replaces the current record by $(x_U,b,\ell+1)$. A constant binding has a read operation returning its stored value. The type checker rejects assignment to a binding of color $C_t$.

For a call with actual argument $(x_U,a,\ell)$ in a live state, the parameter mode determines the binding visible as $y_V$ in the function body. Write $y_{UV}$ for a new local binding introduced on this call from $U$ to $V$. Each invocation creates its own instance of such a binding, including recursive invocations at the same source location.
\begin{description}
    \item[Passing by value.] Initialize a fresh mutable record $(y_{UV},a,0)$. A write through $y_V$ updates this record. The caller's $x_U$ retains its own value and lifetime.
    \item[Passing by reference.] Bind $y_V$ directly to the caller's existing binding $x_U$. Reads and writes through $y_V$ use the current record of $x_U$, including its lifetime state. The call frame records this association for the duration of the call. Releasing the formal binding leaves the lifetime of the caller's object unchanged.
    \item[Constant value binding.] Read $a$ at the call and bind $y_V$ to the pair $(y_{UV},a)$. Each read returns this captured value, including after a subsequent update of $x_U$.
\end{description}

The extractor records whether a formal is a reference and whether it permits writes to the caller's binding. In particular, a constant value parameter owns its copy, whereas a \texttt{const T\&} formal refers to an existing object with writeback disabled. The source interpretation records mutable-reference parameters as well. The current \Lean{} translator rejects functions with such a formal.

On return, the callee first evaluates the return expression, if present, and saves the result for the caller. It then clears its remaining live local slots, including each copy created by passing an argument by value. If the copy $y_{UV}$ last held $(b,\ell')$, its slot undergoes $(b,\ell')\xrightarrow{\mathrm{clear}}\ast$, ending the lifetime of that copy. Releasing the frame discards its constant value bindings and its associations between reference parameters and caller bindings. The caller resumes with the saved result, if any, and sees the current values and lifetime states of its own slots, including any updates made through references. A subsequent call initializes its own by-value copies at state $0$.

\subsection{The Compiler Workflow}
\label{subsec:compiler-workflow}

\begin{figure}[htbp]
    \centering
    \includegraphics[width=\linewidth]{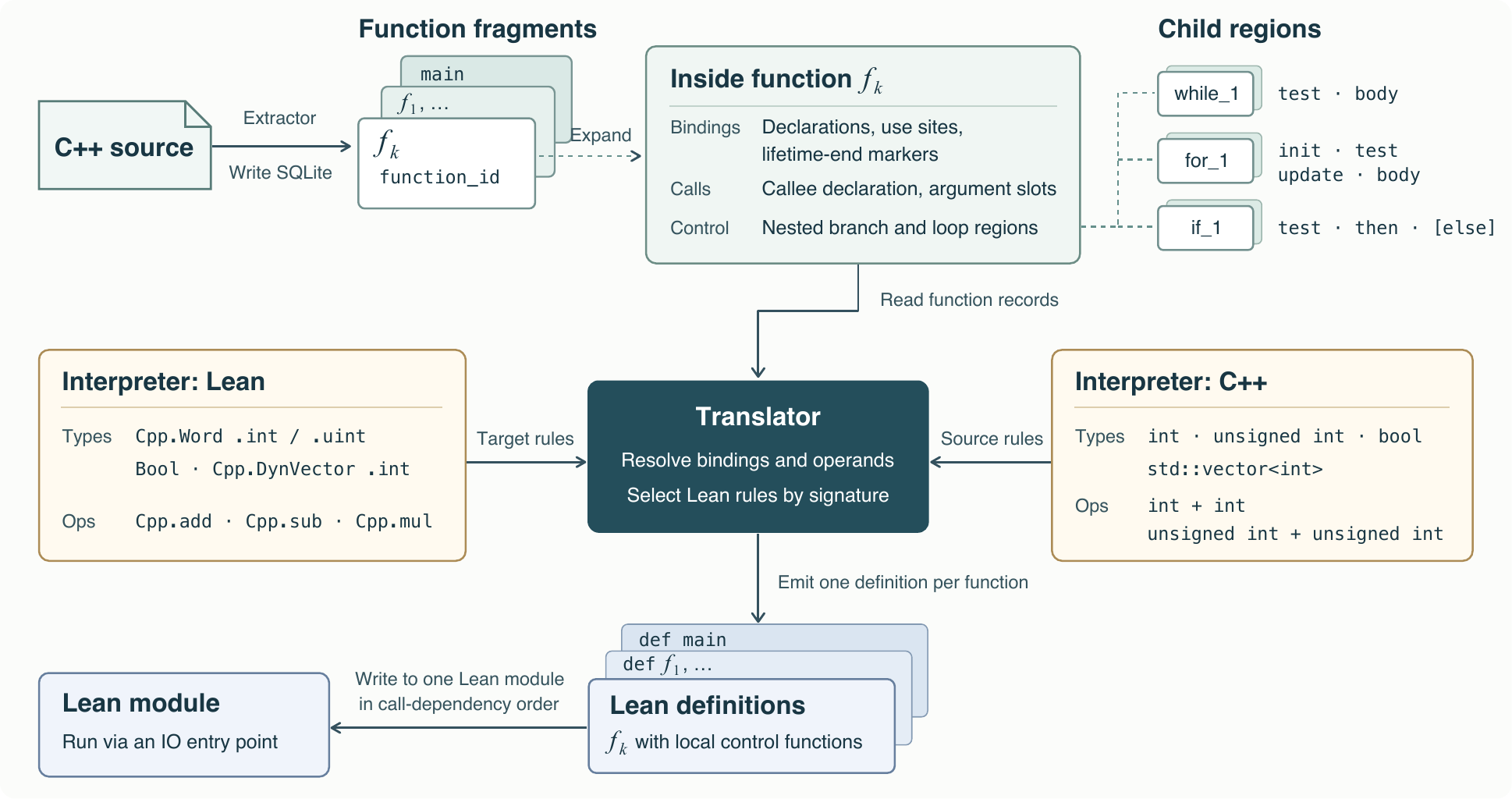}
    \caption{Translation of a function fragment. The top row shows the records for $f_k$ in SQLite. The middle row shows how the translator uses the \Cxx{} and \Lean{} interpretations. The bottom row shows how the resulting definitions form one \Lean{} module.}
    \Description{A function identifier links to bindings, calls, and nested control regions in SQLite. The translator reads these records, resolves operands against the source and target rules, and emits one definition per function. Declaration links determine the order of definitions in the Lean module.}
    \label{fig:compiler-workflow}
\end{figure}

Figure~\ref{fig:compiler-workflow} follows a single function $f_k$ from its \Cxx{} definition to its place in the \Lean{} module. The extractor writes its fragment to SQLite as in Section~\ref{subsec:function-fragments}. Reading the ordered records, the translator resolves variable occurrences to their bindings (Section~\ref{subsec:function-body-translation}) and uses the operator and its operand types to look up the corresponding \Lean{} operation in the interpreter. For example, \texttt{+} on two \texttt{int} operands and \texttt{+} on two \texttt{unsigned int} operands have distinct signatures. In each case the translated operands occupy the same input positions as in the source operation.

The child regions supply the local control functions of Section~\ref{subsec:control-flow-translation}. At each function call, the declaration identifies the target definition and the argument positions determine the parameter bindings of Section~\ref{subsec:parameter-binding}. The translator emits one definition per function and writes them into a single \Lean{} module in call-dependency order, with each callee available to its callers.

Under the hypotheses of Theorem~\ref{thm:faithful-compiler-construction}, the translation works in both directions. With \Cxx{} as $L_1$ and \Lean{} as $L_2$, the theorem gives an isomorphism of semantic algebras $c:A_{L_1}\xrightarrow{\sim}j^*j_!A_{L_1}$. For a \Lean{} definition expressed in this image, the reverse translator reads each operation $j(o)$ as the corresponding \Cxx{} operation $o$ and emits a \Cxx{} source file. Since $c^{-1}$ preserves these operations, a \Cxx{} project can use the generated implementation with the input--output guarantee proved for the \Lean{} algorithm on corresponding inputs. Section~\ref{subsec:application-mlkem} illustrates this use with the polynomial multiplication in ML-KEM.

\section{Applications to Mathematical Algorithms}
\label{sec:compiler-applications}

Following the architecture of Section~\ref{sec:implementation}, we implemented a bidirectional compiler between a subset of \Cxx{}20 and its representation in \Lean{}. On the \Cxx{} side, the implementation supports \texttt{bool}, the signed and unsigned forms of \texttt{short}, \texttt{int}, \texttt{long}, and \texttt{long long}, and functions. It also supports pairs, tuples, fixed-size arrays, bitsets, vectors of integers, and fixed-size Eigen matrices with integer entries.
For the supported fundamental \Cxx{} types and their basic operations, we verify the conditions of Theorem~\ref{thm:faithful-compiler-construction} by showing that translation to \Lean{} preserves values and operation results, including undefined outcomes. Since the case analysis over types and operations is lengthy, we give the full verification in the companion guide on Zenodo~\cite[Sects.~3.4--3.6]{gu-cpp2lean-guide}.

The compiler currently requires \Cxx{} source code to be free of raw pointers. Scientific algorithms that use pointers to traverse arrays can usually be written with indexed access and still perform the same computation. For example, in the Milnor multiplication routine of Section~\ref{subsec:application-milnor}, we read the eight exponents of a monomial by index and retain the original matrix enumeration.

The two applications below arose from practical needs in our own research. First, once we have proved a mathematical algorithm correct, we want to run it efficiently. We believe that AI tools can do much of the work of implementing mathematical algorithms, but the code they generate can be too large to check line by line. Starting from a mathematical specification, we can ask an AI tool to implement the algorithm in \Lean{} and supply a correctness proof for \Lean{} to check. The tool can then rewrite the verified implementation using the types and operations supported by our compiler. Once \Lean{} has checked that the two definitions agree on all admissible inputs, we can compile the rewritten definition to \Cxx{} for efficient execution. For the polynomial multiplication used in ML-KEM, we compiled the verified \Lean{} definition to \Cxx{}. On an Apple M4, the generated code ran \(6.10\) times as fast as the native \Lean{} implementation.

Second, when a mathematical result depends on a large software system, we want to know that the program computes the quantity used in the argument. We therefore used the compiler in the other direction to examine the Milnor multiplication routine used by Lin, Wang, and Xu in their work on the Kervaire invariant one problem in dimension~126~\cite{lin-wang-xu-sseq}. Translating its matrix search to \Lean{} exposed a backtracking step that reads outside the initialized triangle of an array. Repairing this read left the observed products unchanged. The GAIR team's formalization of the Milnor coproduct in \texttt{Kervaire126}~\cite{gair-kervaire126} gives the mathematical specification for this computation.
We present selected results from Sections~5.4 and~6 of the companion guide.

\subsection{Compiling ML-KEM Polynomial Multiplication from \Lean{} to \Cxx{}}
\label{subsec:application-mlkem}

The polynomial multiplication used in ML-KEM takes place in \(R=\mathbb F_{3329}[X]/(X^{256}+1)\). Its inputs and output are vectors of 256 coefficients. For example, \(X^{255}\cdot X=-1\) in \(R\), so the output has constant coefficient \(3328\) and all its other coefficients are zero. We use the number-theoretic transform (NTT) specified in FIPS~203~\cite[Algorithms~9--12]{nist-fips-203}. Microsoft's \texttt{SymCrypt} project provides an executable \Lean{} specification of these operations~\cite{symcrypt-mlkem}. Using Aeneas to translate their Rust implementation into \Lean{}, Ho et al.~\cite{ho-symcrypt-verification} prove that, for any input array of 256 coefficients in \(\{0,\ldots,3328\}\), the forward NTT returns the coefficients prescribed by this specification.

To express the algorithm in the image of the compiler of Theorem~\ref{thm:faithful-compiler-construction}, we take a coefficient vector to be a finite product of 32-bit unsigned words. Each word represents a residue by its unique integer in \(\{0,\ldots,3328\}\). A transform step reads a pair \(u,v\) from the preceding vector and writes \(u+zv\) and \(u-zv\) into the next vector, where \(z\) is the root assigned to that pair. Seven such layers split \(X^{256}+1\) into 128 quadratic factors. Modulo each factor \(X^2-\gamma\), the product of \(a_0+a_1X\) and \(b_0+b_1X\) has coordinates
\[
  \bigl(a_0b_0+\gamma a_1b_1,\;a_0b_1+a_1b_0\bigr).
\]
Applying the inverse transform to these 128 products gives the word vector returned by the \Lean{} definition \texttt{mlkem\_mul}.

For a canonical word vector \(a\), let \(P_a\in R\) be the residue class of the polynomial with those coefficients. For every pair of canonical inputs \(a,b\), we prove in \Lean{} that the word vector \(c=\texttt{mlkem\_mul}(a,b)\) is canonical and that
\[
                         P_c=P_aP_b \quad\text{in }R.
\]
The proof first identifies the word operations with the corresponding field operations. For example, multiplying two canonical words gives a value at most \(3328^2<2^{32}\) before reduction. It then shows that each layer agrees with the sequential updates in our transcription of the specified transform. The forward layers reduce each input modulo the quadratic factors, the displayed formula multiplies the two residues, and the inverse layers recover the product in \(R\).

The same word definition is the input to the reverse translator. It emits a \Cxx{}20 function on \texttt{std::array<unsigned int, 256>}, translating each projection to an indexed read and each bounded fold to a loop. We checked the two native executables on all \(256^2=65{,}536\) ordered pairs of monomials and compared every output coefficient with independent schoolbook multiplication on 67 additional polynomial pairs. On these 67 input pairs, the median time per product over twelve rounds on an Apple M4 was \(57.284\,\mu\mathrm{s}\) for native \Lean{} and \(9.394\,\mu\mathrm{s}\) for generated \Cxx{}. Thus on the same inputs the generated code ran \(6.10\) times as fast. Section~5.4 of the companion artifact guide gives the proof identifiers, test inputs, and timing protocol.

\subsection{Compiling Milnor Multiplication from \Cxx{} to \Lean{}}
\label{subsec:application-milnor}

Lin, Wang, and Xu use the public multiplication operator of \texttt{SSeqCpp} to multiply finite sums of Milnor basis elements in the mod~2 Steenrod algebra~\cite{lin-wang-xu-sseq}. Write \(\mathrm{Sq}(R)\) for the basis element dual to \(\xi^R=\prod_{i=1}^8\xi_i^{r_i}\), of degree \(\sum_{i=1}^8(2^i-1)r_i\). Milnor's coproduct~\cite{milnor-steenrod} gives a specification of the output independent of the program. The coefficient of \(\mathrm{Sq}(T)\) in \(\mathrm{Sq}(R)\mathrm{Sq}(S)\) is the coefficient of \(\xi^R\otimes\xi^S\) in \(\Delta(\xi^T)\), where
\[
                   \Delta(\xi_k)=\sum_{i=0}^{k}\xi_{k-i}^{\,2^i}\otimes\xi_i,
                   \qquad \xi_0=1 .
\]

The public operator calls \texttt{MulMilnorV3} on each ordered pair of summands. It enumerates admissible matrices, packs the resulting monomials into 64-bit words, then sorts the accumulated words and cancels equal pairs. We kept this matrix search and made two changes to its source. First, the original packing routine reads eight consecutive exponents through a pointer. Since our compiler does not currently support pointers, we rewrote the routine to read the same eight exponents by array index. The two versions give identical packed words on all 18,909,086 input monomials through degree~383. The second change concerns the search boundary.

The backtracking array has 81 slots, of which the search initializes only the triangle \(i+j\leq 8\). At \((i,j)=(8,1)\) the original loop tests slot \(9i+j=73\) before checking the boundary. An uninitialized read yields \texttt{none} in the \Lean{} interpretation of the array. We guarded the read with the test \(i+j\leq8\). Inside the triangle the program makes the same test as before. At \((8,1)\) the repaired search returns its accumulator immediately. The \Lean{} proof shows that either Boolean outcome of the original zero test leads to this same accumulator. The zero branch takes one further backtracking step and appends no term.

For example, let \(e_i\) denote the \(i\)th coordinate vector. With a nonunit right factor, the formula gives
\[
       \mathrm{Sq}(2e_1+e_8)\,\mathrm{Sq}(e_1)
       =\mathrm{Sq}(e_2+e_8)+\mathrm{Sq}(3e_1+e_8).
\]
At this input of combined degree~258, \Cxx{} constant evaluation diagnoses the uninitialized read. The original native executable returned the two displayed terms, and the repaired \Lean{} translation returns them as well.

We verified the correctness of the program on all \(1{,}503{,}799{,}211\) ordered pairs of Milnor basis elements with \(\deg R+\deg S\leq180\) by checking that each computed product agrees with the formalized specification above. Section~6 of the companion guide~\cite{gu-cpp2lean-guide} gives the verification records and instructions for reproducing the check.

\section{Related Work and Conclusion}
\label{sec:related-work-conclusion}

\subsection{Related Work}

In the initial algebra semantics of Goguen, Thatcher, Wagner, and Wright~\cite[Sec.~2]{goguen-initial}, the algebra of well-typed ground terms is initial among algebras for the chosen signature. An interpretation of the operation symbols therefore determines a unique homomorphism assigning meanings to the terms. In our setting the value domains of a language and the action of its basic operations are given. We construct $FA_{\mathcal L}$ by adjoining formal expressions and imposing the relations that record the results of computations, including undefined outcomes. The language values embed in this quotient, and evaluating their representatives recovers the original values, as expressed by $\operatorname{res}\circ\iota=\mathrm{id}_{A_{\mathcal L}}$ in Proposition~\ref{prop:semantic-restriction-and-embedding}.

Thatcher, Wagner, and Wright~\cite{thatcher-wagner-wright-compiler} use the initiality of the source syntax in a compiler correctness proof. Once the encoding of source meanings is a homomorphism, there are two homomorphisms from source syntax to target meanings. One compiles a term and interprets the result. The other interprets the term and encodes its meaning. Initiality makes them equal. Mosses~\cite{mosses-constructive} observes that such an encoding can identify distinct source meanings and discusses requiring an injective encoding or specifying a decoding map. Janssen~\cite{janssen-translations} examines the choice between encoding and decoding when a source value has several target representations.

Our compiler starts with a correspondence on types, basic operations, and value generators. Theorem~\ref{thm:faithful-compiler-construction} gives conditions under which this correspondence induces an embedding of the source semantic algebra into the target that preserves both defined and undefined computations. At the level of a faithfully represented program $X$, Proposition~\ref{prop:compiled-program-semantic-equivalence} identifies the algebra recovered by this translation. When the program algebra admits extension by zero along the compiler, restricting the compiled algebra to $\mathcal O_X$ recovers $A_X$ up to canonical isomorphism.

The connection with Lawvere's functorial semantics~\cite{lawvere-functorial}, in its many-sorted form~\cite{benabou-many-sorted}, is one reason we regard algebras over colored operads as a natural framework for program semantics. A colored operad \(\mathcal O\) generates a many-sorted Lawvere theory \(\mathbb L_{\mathcal O}\) with one sort for each color. Restriction induces an equivalence between the set-valued models of \(\mathbb L_{\mathcal O}\) and the set-valued \(\mathcal O\)-algebras~\cite[Sec.~4.6 and Thm.~5.1]{hyland-elements-theories}. In particular, our program algebra \(A_X\) determines a finite-product-preserving functor \(\mathbb L_{\mathcal O_X}\to\mathbf{Set}\). It sends a list of program types to the Cartesian product of their value sets and interprets each program operation by its given action on those values.

In CompCertO, Koenig and Shao~\cite{koenig-shao-compcerto} relate a source component to its compiled code by describing how each interacts with its environment. Their semantics allows correctness proofs for separate components to be composed when the components are linked. Proposition~\ref{prop:program-factorization} recovers the semantic algebra of a single program as the join of the semantic algebras of its function fragments.

Laurent Lafforgue's topos-theoretic view of the semantics of geometric theories~\cite[slides~19--24]{lafforgue-topos-semantics} inspired us to use a factorization algebra to describe a program. We take the program \(X\) itself as the underlying space and assign to each open fragment \(U\) its semantic algebra \(\mathcal F_X(U)\). For the cover by function fragments, Proposition~\ref{prop:program-factorization} recovers \(A_X\) as the join of these local algebras under their embeddings in \(A_X\).

\subsection{Conclusion}

We have described the semantics of a program by an algebra over a colored operad, independently of the language in which the program is expressed. Semantic restriction relates formal expressions to their values, including undefined outcomes, and a semantically faithful translation preserves the distinction between values of the same type. The same algebra can be viewed as a model of the associated many-sorted Lawvere theory. Working directly with the colored operad and its algebra makes the conditions for preserving computations easier to check. We check that the translation respects the defining relations and commutes with each basic operation. Compatibility with every composite then follows by operadic substitution. For a program faithfully represented in the source language, Proposition~\ref{prop:compiled-program-semantic-equivalence} shows that, when the required extension by zero exists, restricting the compiled algebra to the program's own types and operations recovers its original semantic algebra up to canonical isomorphism.

Once a compiler has been verified on fundamental types and operations, which parts of that verification must be repeated for a new type or operation? For a product, the translations of the coordinates determine the translation of the tuple. For recursion, it suffices to check that the translations preserve the initial and step operations. Theorem~\ref{thm:compiler-extension-products-recursion} shows that any finite sequence of these extensions preserves semantic faithfulness. Reconstruction from function fragments gives a corresponding organization for the implementation. The \Cxx{}-to-\Lean{} compiler translates each function body once, and each call refers to the resulting definition.

For the polynomial multiplication used in ML-KEM, we expressed the number-theoretic transform by finite products and bounded recursion. This gave us a single \Lean{} definition on which to prove the ring identity for all canonical inputs and from which to generate an executable \Cxx{} implementation. For the Milnor multiplication used in the dimension~126 Kervaire invariant one computation~\cite{lin-wang-xu-sseq}, we started from an existing \Cxx{} implementation. Its \Lean{} interpretation made an uninitialized read visible during backtracking. After repairing the boundary test, we verified that the program computes the product prescribed by Milnor's coproduct on every ordered basis pair of combined degree at most~180.

\begin{acks}
The author is deeply grateful to Guozhen Wang for suggesting this problem and for his many years of guidance and support as an advisor. Discussions with him inspired both the use of finite products and recursion to extend the compiler in Section~\ref{subsec:compiler-extension} and the treatment of branches and loops as functions in Section~\ref{subsec:control-flow-translation}.
The author thanks Weinan Lin for his many suggestions on an early draft and for the idea of compiling verified \Lean{} programs to \Cxx{} to improve their performance.
Thanks are also due to Zhengping Gui, whose course introduced the author to factorization algebras and their applications in mathematics and physics, and to Shuhan Zheng for discussions on this subject.
The author is grateful to Laurent Lafforgue for an in-person conversation in which he explained how the Langlands program led him to describe languages using Grothendieck toposes. His perspective inspired the author's use of factorization algebras in program semantics.
\end{acks}

\bibliographystyle{ACM-Reference-Format}
\bibliography{references}

\end{document}